\documentclass{amsart}
\usepackage{cmap}
\usepackage[T1]{fontenc}
\usepackage[utf8]{inputenc}
\usepackage{lmodern}

\usepackage{amsmath}
\usepackage{amssymb}
\usepackage{amsthm}

\usepackage{enumitem}
\usepackage{xcolor}

\usepackage{mathrsfs}
\usepackage{caption}
\usepackage{subcaption}
\usepackage{tikz-cd}
\usetikzlibrary{arrows}
\usetikzlibrary{arrows.meta}
\usetikzlibrary{decorations.markings}
\usetikzlibrary{calc, 3d}
\usepackage{hyperref}

\hypersetup{
	pdftitle={Deleted square of a cone over a graph},
	pdfauthor={Karolak Andrzej T.},
	pdfkeywords={configuration space, deleted product, deleted square, deleted join, graph, cone}
}
\newtheorem{Thm}{Theorem}[section]

\newtheorem{Cor}[Thm]{Corollary}
\newtheorem{Lem}[Thm]{Lemma}
\newtheorem{Prop}[Thm]{Proposition}
\newtheorem{Conj}[Thm]{Conjecture}
\theoremstyle{definition}
\newtheorem{Rem}[Thm]{Remark}

\DeclareMathOperator{\Conf}{Conf}

\DeclareMathOperator{\UConf}{UConf}
\DeclareMathOperator{\supp}{supp}
\DeclareMathOperator{\coker}{coker}
\DeclareMathOperator{\st}{st}
\DeclareMathOperator{\id}{id}
\DeclareMathOperator{\Cn}{Cn}
\DeclareMathOperator{\Lk}{Lk}
\DeclareMathOperator{\adv}{AV}
\DeclareMathOperator{\ade}{AE}
\DeclareMathOperator{\gmc}{GMC}
\DeclareMathOperator{\colim}{colim}
\DeclareMathOperator{\hoco}{hocolim}
\newcommand{\sdp}{\mathbin{\times_{\triangle}}}
\newcommand{\sdj}{\mathbin{*_{\triangle}}}

\newcommand{\G}{\Gamma}
\newcommand{\Z}{\mathbb{Z}}
\newcommand\quot[2]{
	\mathchoice
	{
		\text{\raise1ex\hbox{$#1$}\Big/\lower1ex\hbox{$#2$}}%
	}
	{
		#1\,/#2
	}
	{
		#1\,/#2
	}
	{
		#1\,/#2
	}
}

\begin{document}
	\title{Deleted square of a cone over a graph}
	\author{Andrzej T. Karolak}
	\address{\begin{tabular}[t]{@{}l@{}}
			Instytut Matematyczny\\ Uniwersytet Wroc\l{}awski\\
			pl. Grunwaldzki 2\\ 50-384 Wroc\l{}aw\\ Poland
	\end{tabular}}
	\email{andrzej.karolak@math.uni.wroc.pl}
	\begin{abstract}
		The deleted square of a topological space is a subspace of the Cartesian product of this space with itself, obtained by removing the diagonal. Conjecturally, the deleted square of a cone over a graph has the homotopy type of a wedge sum of spheres of dimensions at most 3. We propose a geometric construction that is a key motivation for the conjecture and provides a major step towards its proof, while yielding direct results on deleted squares of certain families of graphs.
	\end{abstract}
	
	\maketitle
	\section{Introduction}
	In this paper, we study the deleted square of a cone over a graph.
	The \emph{deleted square} of a topological space $X$ is defined as the following subspace of the Cartesian product of $X$ with itself:
	\[X^2_{\circ}=\{(x_1,x_2)\in X^2\mid x_1\neq x_2\}.\]
	
	The deleted square is an instance of an $n$-fold \emph{deleted product} of a space $X$, also known as a \emph{residual space}, defined as the following subspace of the Cartesian product $X^n$:
	\[X^n_{\circ}=R_n(X)=X^n\setminus \Delta_n(X),\] where $\Delta_n(X)\colon X\to X^n$ is the embedding of the diagonal: $\Delta_n(x)=(x,\dots,x)$.
	
	It is also a special case of the \emph{configuration space} of $n$ particles on $X$, which is defined as the following subset of $X^n$:
	\[\Conf_n(X)=\{(x_1,\dots,x_n)\in X^n\mid i\neq j \implies x_i\neq x_j\},
	\]
	which is often called an \emph{ordered} configuration space as opposed to the \emph{unordered} one, often denoted as $\UConf_n(X)$, which is a quotient of $\Conf_n(X)$ by the natural action of the symmetric group permuting the coordinates.
	
	One of the major milestones in the research on configuration spaces on graphs was achieved by Abrams \cite{abrams2000configuration} and \'Swi\k{a}tkowski \cite{Swiatkowski}. In these independent papers, cube complexes that are deformation retracts of $\Conf_n(\G)$ and of $\UConf_n(\G)$ were described. Another result on the topology of configuration spaces was presented in a paper by Gal \cite{Gal}, in which the author computed the Euler characteristic of the configuration space on an arbitrary simplicial complex. Surprisingly, it was still a challenge to compute the Betti numbers of the configuration space on a graph. The case of two particles was ultimately settled by Ko and Park \cite{ko2012characteristics}, who computed $H_1(\G^2_{\circ})$ using discrete Morse theory. Prior to their work, Farber and Barnett~\cite{barnett2009topology}, along with Farber and Hanbury~\cite{farber2010topology}, investigated the two-particle case in a two-part series of papers, describing its homology for certain families of graphs.
	
	In all these papers, only algebraic invariants are computed explicitly. It turns out that in the case of the deleted square of a cone over a graph, more can be said about the homotopy type itself.
	
	We propose a construction that describes how the deleted square of a cone over a graph changes when the underlying graph is modified. It was inspired by the graph enlarging construction from \cite{farber2010topology}.
	
	Throughout this paper, every graph will be tacitly assumed to be finite.
	
	Every connected graph can be built step by step, starting from a single vertex, by using two operations:
	\begin{itemize}
		\item \emph{Adding a new vertex as a leaf}: for a graph $\G$ and its vertex $v$, we define a new graph $\adv(\G,v)$, which is obtained by identifying one endpoint of $[0,1]$ with $v$.
		\item \emph{Adding a new edge between existing vertices}: for a graph $\G$ and its distinct vertices $v_1$ and $v_2$, we define a new graph $\ade(\G,v_1,v_2)$, which is obtained by identifying the endpoints of $[0,1]$ with $v_1$ and $v_2$, respectively. 
	\end{itemize}
	
	The effect of these operations on the homotopy type of the deleted square of a cone over a graph is expressed by the theorem:
	\begin{Thm}
		For a graph $\G$, and its vertices $v$, $v_1$, and $v_2$, where $v_1\neq v_2$:
		\begin{enumerate}[label=\alph*)]
			\item [a)] $\Cn(\adv(\G,v))^2_{\circ}$ is homotopy equivalent to the homotopy colimit of the following diagram:
			\[\begin{tikzcd}
				* & X_1 \arrow{l}\arrow{r} & (\Cn\G)^2_{\circ} & X_2 \arrow{l}\arrow{r} & *.
			\end{tikzcd}\]
			\item [b)] $\Cn(\ade(\G,v_1,v_2))^2_{\circ}$ is homotopy equivalent to the homotopy colimit of the following diagram:
			\[\begin{tikzcd}
				* & X_3 \arrow{l}\arrow{r} & (\Cn\G)^2_{\circ} & X_4 \arrow{l}\arrow{r} & *.
			\end{tikzcd}\]
		\end{enumerate}
		where $*$ denotes the one-point space, $\Cn\G$ denotes the cone over $\G$, and $X_i$'s are certain subspaces of $(\Cn\G)^2_{\circ}$ of following homotopy types:
		\[X_1\simeq X_2 \simeq \bigvee_{i=1}^{\nu(v)-1}S^1,\]
		and
		\[X_3 \simeq X_4 \simeq \left(\bigvee_{i=1}^{b_0-1}S^1\right)\vee\left(\bigvee_{i=1}^{b_1}S^2\right),\]
		where $\nu(v)$ denotes the valence or degree of the vertex $v$ in $\G$, which is the number of vertices of $\G$ adjacent to $v$, and $b_0$ and $b_1$ are the zeroth and the first Betti numbers of $\G\setminus\{v_1,v_2\}$.
	\end{Thm}
	To state this more intuitively, obtaining the deleted square of a cone over a graph $\G$ with an added leaf or edge requires taking separate cones over two particular subspaces of $(\Cn\G)^2_{\circ}$.
	
	This construction led us to the following conjecture:
	\begin{Conj}\label{conj}
		$(\Cn\G)^2_{\circ}$ is homotopy equivalent to a bouquet of 0-, 1-, 2-, and 3-dimensional spheres.
	\end{Conj}
	It is known that
	$(\Cn *)^2_{\circ}\simeq S^0$ and $(\Cn [0,1])^2_{\circ}\simeq(\Cn S^1)^2_{\circ}\simeq S^1$.
	To simplify the assumptions throughout the paper, we call a graph homeomorphic to a point, a line segment, or a circle \emph{degenerate}, and all the other connected graphs \emph{nondegenerate}. Equivalently, a graph is nondegenerate if it is connected and has a vertex of valence at least 3.
	
	In Proposition~\ref{discon}, we discuss disconnected graphs, which allows us to restrict our attention to connected graphs.
	
	These distinctions are essential, as the conjecture can be restated as follows:
	\begin{Conj}\label{conj2}
		If $\G$ is a nondegenerate graph, then $(\Cn\G)^2_{\circ}$ is homotopy equivalent to a bouquet of 2- and 3-dimensional spheres.
	\end{Conj}
	
	This conjecture provides more insight into the facts proven in \cite{barnett2009topology} and \cite{farber2010topology}, and some properties of the deleted squares of graphs can be deduced from the study of the deleted squares of their cones. For example, Proposition~2.3 in \cite{barnett2009topology} states:
	\begin{Prop}[Farber, Barnett]
		Let $\G$ be a nondegenerate graph. Then:
		\begin{enumerate}[label=\roman*)]
			\item $H_2(\Conf_2(\G)) \cong \ker(I_{\G})$,
			\item $H_1(\Conf_2(\G)) \cong \coker(I_{\G}) \oplus H_1(\G) \oplus H_1(\G)$,
		\end{enumerate}
		where the map \[I_{\G}\colon H_1(\G) \otimes H_1(\G)\to H_2(\G^2,\Conf_2(\G))\] appears in the long exact sequence of homology of the pair $(\G^2,\Conf_2(\G))$.
	\end{Prop}
	
	A similar description is obtained by examining the Mayer--Vietoris sequence of the configuration space of a cone over a graph.
	\begin{Prop}[Copeland]\label{propA}
		Let $\G$ be a nondegenerate graph. Then:
		\begin{enumerate}[label=\roman*)]
			\item $H_2(\Conf_2(\G)) \cong H_3(\Conf_2(\Cn\G))$.
			\item $H_1(\Conf_2(\G)) \cong H_2(\Conf_2(\Cn\G)) \oplus H_1(\G) \oplus H_1(\G)$.
		\end{enumerate}
	\end{Prop}
	
	This proposition is a special case of a result of Copeland \cite{copeland1965deleted}, which shows a relation between the cohomology of $X$, of the deleted square of $X$, and of the deleted square of a cone over $X$ for any triangulated space $X$. We present a proof adapted to the presented case in Section~3.
	
	In \cite{farber2010topology}, Farber and Hanbury point to a trichotomy present in the properties of the map
	\[\alpha_*\colon H_1(\G^2_{\circ})\to H_1(\G^2)\]
	induced by the inclusion $\alpha\colon\G^2_{\circ}\to\G^2$:
	\begin{enumerate}

		\item In the case of trees, $\alpha_*$ has a large kernel.
		\item In the case of planar graphs with all vertices of valence at least $3$ and some additional conditions (see \cite{barnett2009topology}), $\alpha_*$ has a kernel isomorphic to $\Z$.
		\item For some graphs, $\alpha_*$ is an isomorphism.
	\end{enumerate}
	The authors called the graphs satisfying (3) \emph{mature}.
	
	In almost all cases, $\alpha_*$ is an epimorphism, so this is in fact a statement about the kernel of $\alpha_*$, which coincides with $\coker I_{\G}$ and (as Proposition~\ref{propA} indicates) with $H_2(\Conf_2(\Cn\G))$. Utilizing this fact, we propose a similar trichotomy based on the conjecture.
	For a nondegenerate graph $\G$, $(\Cn\G)^2_{\circ}$ has the homotopy type of:
	\begin{enumerate}
		\item a bouquet of 2-spheres,
		\item a bouquet of 2- and 3-spheres,
		\item a bouquet of 3-spheres.
	\end{enumerate}
	We call a graph $\G$ for which (1) holds \emph{unripe}, and a graph for which (3) holds \emph{ripe}.
	
	Ripe graphs are obviously mature, and the two notions would be equivalent under the assumption that Conjecture~\ref{conj} holds. 
	
	We show that trees and complete bipartite graphs of the form $K_{2,m}$ with $m > 2$ are unripe and that certain theorems on mature graphs from \cite{farber2010topology} can be translated to ripe graphs. In particular, complete graphs $K_n$ with $n \geq 5$ and complete bipartite graphs $K_{k,m}$ with $k,m \geq 3$ are ripe.
	
	The paper is organized as follows. In Section~2, we establish our tools and terminology by recalling basic constructions such as the cone, introducing the notion of the deleted square, and developing the necessary framework for homotopy colimits. In Section~3, we present basic facts about the deleted square of a cone over a graph, including a proof of Proposition~\ref{propA}. To give an example, we compute the homotopy type of $(K_{n,m})^2_{\circ}$. In Section~4, we build on the idea from \cite{farber2010topology} of enlarging a graph by attaching an edge to it, and prove the main theorem. In Section~5, we utilize the main theorem to prove facts about the homotopy type of deleted squares of cones over certain graphs, and draw a comparison between them and the results of \cite{barnett2009topology} and \cite{farber2010topology}. Section~6 is an appendix consisting of additional proofs and theorems included for completeness.
	
	\section {Preliminaries}
	\subsection{Basic notions} We recall some basic notions used throughout the paper mainly to establish notation. Most of the constructions in this section could be defined generally for regular cell complexes, but we restrict our attention to simplicial complexes.
	
	Whenever we talk about a simplicial complex $X$, we consider it primarily as a topological space, but we also recognize its combinatorial structure, and we abuse notation relying on the context. For example, if we write $\sigma\in X$, we consider $\sigma$ as a combinatorial simplex of $X$, so as a subset of the set of vertices of $X$, but if we write $\sigma\subset X$, we consider this simplex as a topological subspace. In other words, we identify a simplicial complex with its geometric realization.
	
	We consistently use the symbol $\cong$ for homeomorphism and $\simeq$ for homotopy equivalence. 
	
	We define a \emph{cone} over a topological space $X$ as the quotient:
	\[\Cn X=\quot{X\times[0,1]}{\sim},\]
	where the relation is given by $(x,1)\sim(x',1)$ for all $x,x'\in X$.
	
	For a point $(x,t)$ in $X\times [0,1]$, we denote its class by $[x,t]$, and for any $x\in X$, we call $[x,1]$ the \emph{apex} of the cone. We identify $X$ with the image of $X\times\{0\}$ under the quotient map and call it the \emph{base} of the cone.
	
	If $X$ is a simplicial complex, then $\Cn X$ has a natural structure of a simplicial complex that can be defined via the combinatorial join $X*\Delta^0$, where $\Delta^0$ is the $0$-dimensional simplex (a vertex). In general, for two simplicial complexes $X$ and $Y$, we define their combinatorial join:
	\[X*Y=\{(\sigma\times\{-1\})\cup(\tau\times\{1\})\mid\sigma\in X,\tau\in Y\}.\]
	Taking products with $\{-1\}$ and $\{1\}$ ensures that $X$ and $Y$ are treated as disjoint complexes, even if they share simplices (for instance, when they are subcomplexes of the same complex).
	
	We also define the topological join (we abuse notation without proving that it is homeomorphic to the geometric realization of the combinatorial join):
	\[X*Y=\quot{X\times Y\times[-1,1]}{\sim},\]
	where $\sim$ is given by $(x,y,-1)\sim(x,y',-1)$ and $(x,y,1)\sim(x',y,1)$ for any $x,x'\in X$ and $y,y'\in Y$.
	The choice of $-1$ and $1$ will be convenient for an argument later. We denote points of the join by $[x,y,t]$ and identify $X$ with the image of $X\times Y\times\{-1\}$ and $Y$ with the image of $X\times Y\times\{1\}$ under the quotient map.
	
	The \emph{suspension} $SX$ of a space $X$ is the join of $X$ with the discrete space consisting of two points.
	
	For a simplicial complex $X$, we introduce the following space:
	\[X^{2}_{\triangle}=\{(x,y)\in X^2\mid\supp_X(x)\cap\supp_X(y)=\emptyset\},\]
	where $\supp_X(x)$ denotes the smallest closed simplex of $X$ containing $x$. It is called the \emph{simplicial deleted square} of $X$. It is the subcomplex of $X^2$ consisting of those closed cells that do not intersect the diagonal. More generally, we define
	\[X\sdp Y=\{(x,y)\in X\times Y\mid\supp_X(x)\cap\supp_Y(y)=\emptyset\},\]
	for complexes $X$ and $Y$. In particular, $X^2_{\triangle}=X\sdp X$.
	
	This definition is primarily intended for intersecting subcomplexes of a given complex. However, we intentionally adopt this form to include the case of disjoint complexes $X$ and $Y$, for which $X \sdp Y = X \times Y$.
	
	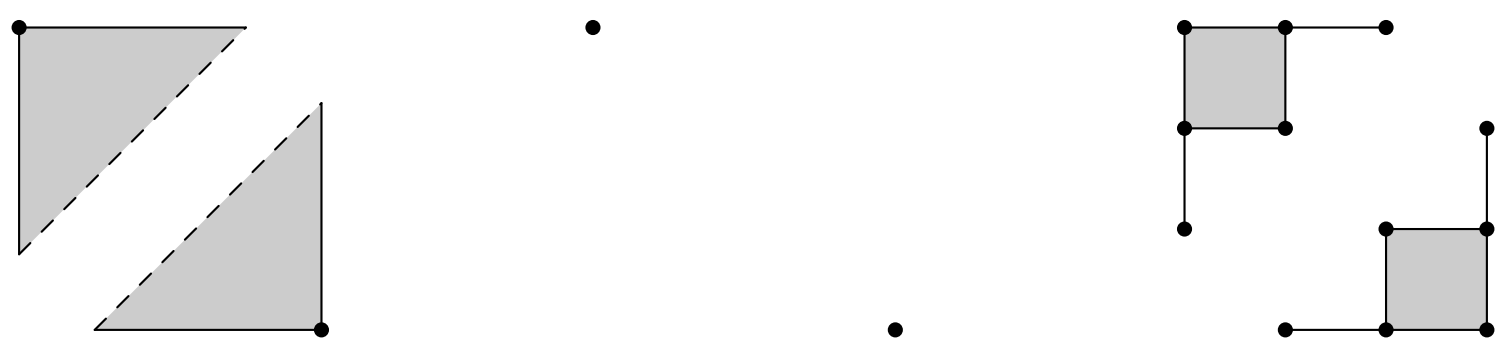
\begin{figure}
		\centering
		\begin{subfigure}{.3\textwidth}
			\centering
			\begin{tikzpicture} [line cap=round,line join=round,x=1.0cm,y=1.0cm]
				\clip (-1.7,-1.5) rectangle (1.7,1.5);
				\draw (-1,-0.5)--(-1,1)--(0.5,1);
				\draw [dashed] (-1,-0.5)--(0.5,1);
				\fill [opacity=0.2] (-1,-0.5)--(-1,1)--(0.5,1)--cycle;
				\draw (-0.5,-1)--(1,-1)--(1,0.5);
				\draw [dashed] (-0.5,-1)--(1,0.5);
				\fill [opacity=0.2] (-0.5,-1)--(1,-1)--(1,0.5)--cycle;
				\fill (-1,1) circle[radius=0.05];
				\fill (1,-1) circle[radius=0.05];
			\end{tikzpicture}
			\caption{$(\Delta^1)^2_{\circ}$}
			\label{fig1a}
		\end{subfigure}%
		\begin{subfigure}{.3\textwidth}
			\centering
			\begin{tikzpicture} [line cap=round,line join=round,x=1.0cm,y=1.0cm]
				\clip (-1.7,-1.5) rectangle (1.7,1.5);
				\fill (-1,1) circle[radius=0.05];
				\fill (1,-1) circle[radius=0.05];
			\end{tikzpicture}
			\caption{$(\Delta^1)^2_{\triangle}$}
			\label{fig1b}
		\end{subfigure}
		\begin{subfigure}{.3\textwidth}
			\centering
			\begin{tikzpicture} [line cap=round,line join=round,x=1.0cm,y=1.0cm]
				\clip (-1.7,-1.5) rectangle (1.7,1.5);
				\draw (-1,-0.333)--(-1,1)--(0.333,1);
				\draw (-1,0.333)--(-0.333,0.333)--(-0.333,1);
				\draw (-0.333,-1)--(1,-1)--(1,0.333);
				\draw (0.333,-1)--(0.333,-0.333)--(1,-0.333);
				\foreach \p in {
					(-1,-0.333), (-1,1), (0.333,1),
					(-1,0.333), (-0.333,0.333), (-0.333,1),
					(-0.333,-1), (1,-1), (1,0.333),
					(0.333,-1), (0.333,-0.333), (1,-0.333)
				} \fill \p circle[radius=0.05];
				\fill[opacity=0.2] (-1,0.333) rectangle (-0.333,1);
				\fill[opacity=0.2] (0.333,-1) rectangle (1,-0.333);
			\end{tikzpicture}
			\caption{$Sd^2_{\triangle}$}
			\label{fig1c}
		\end{subfigure}
		\caption{(A) The deleted square and (B) the simplicial deleted square of a $1$-simplex $\Delta^1$; (C) the simplicial deleted square of a subdivision $Sd$ of $\Delta^1$ into three edges.}
		\label{fig:line}
	\end{figure}
	
	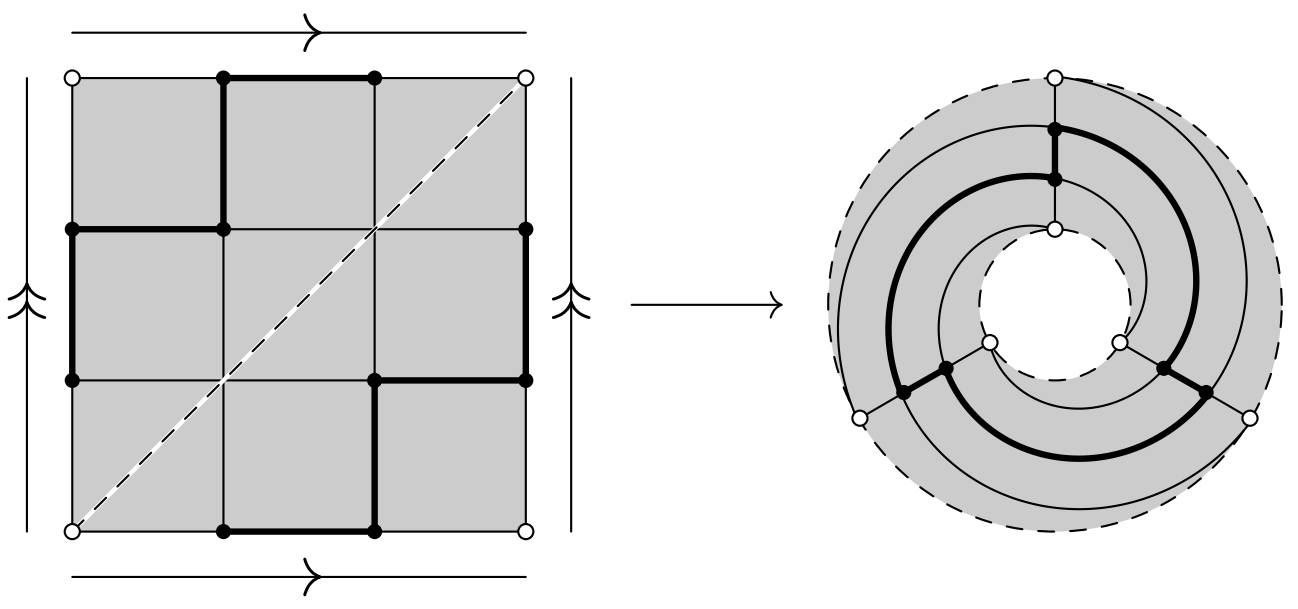
\begin{figure}
		\centering
		\makebox[\textwidth][c]{%
			\begin{tikzpicture} [line cap=round,line join=round,x=1.0cm,y=1.0cm, arrowA/.style={
					postaction={decorate},
					decoration={
						markings,
						mark=at position 0.55 with {\arrow[ scale=1.4]{>}}
					}
				},
				arrowB/.style={
					postaction={decorate},
					decoration={
						markings,
						mark=at position 0.55 with {\arrow[ scale=1.4]{>>}}
					}
				}
				]
				\clip (-1,-0.5) rectangle (9,3.5);
				\fill [opacity=0.2] (0,0) rectangle (3,3);
				\draw (0,0) rectangle (3,3);
				\draw (0,1)--(3,1);
				\draw (0,2)--(3,2);
				\draw (1,0)--(1,3);
				\draw (2,0)--(2,3);
				\draw [thick, white] (0,0)--(3,3);
				\draw [dashed] (0,0)--(3,3);
				
				\draw [arrowA] (0,-0.3) -- (3,-0.3);
				\draw [arrowA] (0,3.3) -- (3,3.3);
				
				\draw [arrowB] (-0.3,0) -- (-0.3,3);
				\draw [arrowB] (3.3,0) -- (3.3,3);
				
				\foreach \x/\y in {0/0, 0/3, 3/0, 3/3} {
					\fill [color=white] (\x,\y) circle [radius=0.05];
					\draw (\x,\y) circle [radius=0.05];
				}
				
				\draw [very thick] (0,1)--(0,2)--(1,2)--(1,3)--(2,3);
				\draw [very thick] (1,0)--(2,0)--(2,1)--(3,1)--(3,2);
				\fill (0,1) circle [radius=0.05];
				\fill (0,2) circle [radius=0.05];
				\fill (1,2) circle [radius=0.05];
				\fill (1,3) circle [radius=0.05];
				\fill (2,3) circle [radius=0.05];
				\fill (1,0) circle [radius=0.05];
				\fill (2,0) circle [radius=0.05];
				\fill (2,1) circle [radius=0.05];
				\fill (3,1) circle [radius=0.05];
				\fill (3,2) circle [radius=0.05];
				
				\draw [->] (3.7,1.5)--(4.7,1.5);
				
				\fill[opacity=0.2, even odd rule] (6.5,1.5) circle (1.5) circle (0.5);
				\draw [dashed] (6.5,1.5) circle [radius=1.5];
				\draw [dashed] (6.5,1.5) circle [radius=0.5];
				\draw (6.5,2)--(6.5,3);
				\draw (6.93,1.25)--(7.79,0.75);
				\draw (6.07,1.25)--(5.21,0.75);

				\draw[domain=90:450, samples=300, shift={(6.5,1.5)}, variable=\t]
				plot (\t:{0.0028*(\t-90)+0.5});
				
				\draw[domain=330:690, samples=300, shift={(6.5,1.5)}, variable=\t]
				plot (\t:{0.0028*(\t-330)+0.5});
				
				\draw[domain=570:930, samples=300, shift={(6.5,1.5)}, variable=\t]
				plot (\t:{0.0028*(\t-570)+0.5});
				
				\fill [color=white] (6.93,1.25) circle [radius=0.05];
				\draw (6.93,1.25) circle [radius=0.05];
				\fill [color=white] (7.79,0.75) circle [radius=0.05];
				\draw (7.79,0.75) circle [radius=0.05];
				
				\fill [color=white] (6.07,1.25) circle [radius=0.05];
				\draw (6.07,1.25) circle [radius=0.05];
				\fill [color=white] (5.21,0.75) circle [radius=0.05];
				\draw (5.21,0.75) circle [radius=0.05];
				
				\fill [color=white] (6.5,2) circle [radius=0.05];
				\draw (6.5,2) circle [radius=0.05];
				\fill [color=white] (6.5,3) circle [radius=0.05];
				\draw (6.5,3) circle [radius=0.05];
				
				\fill (6.5,2.33) circle [radius=0.05];
				\fill (6.5,2.66) circle [radius=0.05];
				\draw [very thick] (6.5,2.33)--(6.5,2.66);
				
				\fill (7.22,1.08) circle [radius=0.05];
				\fill (7.5,0.92) circle [radius=0.05];
				\draw [very thick] (7.22,1.08)--(7.5,0.92);
				
				\fill (5.5,0.92) circle [radius=0.05];
				\fill (5.78,1.08) circle [radius=0.05];
				\draw [very thick] (5.5,0.92)--(5.78,1.08);

				\draw[very thick, domain=210:330, samples=300, shift={(6.5,1.5)}, variable=\t]
				plot (\t:{0.0028*(\t-90)+0.5});
				
				\draw[very thick, domain=90:210, samples=300, shift={(6.5,1.5)}, variable=\t]
				plot (\t:{0.0028*(\t+30)+0.5});
				
				\draw[very thick, domain=330:450, samples=300, shift={(6.5,1.5)}, variable=\t]
				plot (\t:{0.0028*(\t-210)+0.5});
				
		\end{tikzpicture}}
		
		\caption{The deleted square of a hollow triangle $\partial \Delta^2$ with its simplicial deleted square indicated by thick lines.}
		\label{fig:circle}
	\end{figure}
	
	The following result is already known (see \cite{abrams2000configuration} or \cite{hu1960isotopy} for instance), but we give its proof for completeness in the appendix.
	\begin{Lem}\label{lemret}
		$X^2_{\triangle}$ is a strong deformation retract of $X^2_{\circ}$.
	\end{Lem}
	\begin{Rem}
		Notice that the topology of $X^2_{\triangle}$ depends on the combinatorial structure of $X$, but its homotopy type does not, as long as the change in combinatorial structure does not affect the topology of $X$. Indeed, if $X$ is homeomorphic to $X'$, then
		\[X^2_{\triangle}\simeq X^2_{\circ}\cong (X')^2_{\circ}\simeq (X')^2_{\triangle}.\]
		In particular, for any subdivision $X'$ of a complex $X$, we have $X^2_{\triangle}\simeq (X')^2_{\triangle}$.
	\end{Rem}
	\begin{Rem}
		$X^2_{\triangle}$ is clearly not a simplicial complex, but it is still a regular cell complex. We are not concerned with its structure as such. For our purposes, we only need its homotopy-theoretic properties as a CW-complex, that is, the fact that an embedding of a subcomplex in a CW-complex is a cofibration (see \cite[Chapter 7]{kozlov2008combinatorial}), although we don't use this fact explicitly.
	\end{Rem}
	
	\subsection{Homotopy colimit}
	We now take a slight detour to set up an elementary framework for the homotopy-theoretic constructions used in the proof of the main theorem. Although our diagrams are simple enough that all relevant homotopy colimits reduce to generalized mapping cylinders, we rely on several results whose precise statements require the language of homotopy colimits. In particular, we need a homotopy-commutative variant of the standard homotopy-invariance lemma (Lemma~\ref{simhoco}), for which we include an elementary proof in the appendix for completeness. Fixing this notation and terminology avoids appealing to higher category theory and reduces technicalities in the main proof. Following \cite[Section~4.G]{hatcher2005algebraic} and \cite[Chapter~15]{kozlov2008combinatorial}, we recall and adapt the required definitions below.
	
	A \emph{diagram of spaces} consists of a directed graph $\G$, a topological space $X_v$ for each vertex $v$ of $\G$, and a continuous map $\varphi_e\colon X_v\to X_w$ for each edge $e=(v,w)$ of $\G$. We denote this system by $X_\G=(X_v,\varphi_e)_{\G}$. We only consider finite connected graphs $\G$ in this context.
	
	This definition comes from \cite{hatcher2005algebraic}. It is related to the notion of an \emph{arrow picture} in \cite{kozlov2008combinatorial}. The notion of a diagram of spaces is defined in \cite{kozlov2008combinatorial}
	over triangulated spaces (trisps), where higher-dimensional cells indicate commutative parts of the diagram. 
	
	For a diagram of spaces $X_{\G}=(X_v,\varphi_e)_{\G}$, we define its \emph{generalized mapping cylinder} as the following quotient space
	\[\gmc(X_{\G})=\quot{\left(\bigsqcup_{v\in V(\G)}X_v \sqcup \bigsqcup_{e\in E(\G)}X_{0_e}\times [0,1]_e\right)}{\sim},\]
	where $[0,1]_e$ is a copy of the unit interval associated with the edge $e$. To emphasize this, we denote a point of $[0,1]_e$ associated with $t\in[0,1]$ by $t_e$. We identify the initial vertex of $e$ with $0_e$ and the terminal one with $1_e$, so that we can write $\varphi_e\colon X_{0_e}\to X_{1_e}$, where $X_{0_e}=X_v$ and $X_{1_e}=X_w$ if $e=(v,w)$.
	
	The equivalence relation $\sim$ is induced by $(x,0_e)\sim x$ and $(x,1_e)\sim \varphi_e(x)$ for each edge $e$ and each $x\in X_{0_e}$.
	
	We also recall that the \emph{colimit} (\emph{amalgamation} in \cite{hatcher2005algebraic}) of $X_{\G}$ is defined as the quotient
	\[\colim (X_{\G})=\quot{\bigsqcup_{v\in V(\G)}X_v}{\sim},\]
	where $\sim$ is induced by $x\sim\varphi_e(x)$ for each edge $e$ and each $x\in X_{0_e}$.
	
	We call a vertex $s$ of a directed graph $\G$ a \emph{source} if there is no edge of $\G$ of the form $(v,s)$. In other words, all the edges adjacent to $s$ point outward from it. Similarly, we call a vertex $d$ of $\G$ a \emph{drain} if there is no edge of $\G$ of the form $(d,v)$.
	We call a diagram of spaces over $\G$ a \emph{source--drain diagram} if every vertex of $\G$ is either a source or a drain.
	
	We only consider source--drain diagrams, and in this case, the homotopy colimit of a source--drain diagram of spaces is homeomorphic to the generalized mapping cylinder of this diagram, which follows directly from the definition of a homotopy colimit (see \cite[page~262]{kozlov2008combinatorial}). Thus, whenever we refer to a homotopy colimit of a source--drain diagram, we mean either its generalized mapping cylinder or a space homotopy equivalent to it. To make a distinction, $\gmc$ will be referred to as \emph{the} homotopy colimit, and the others as \emph{a} homotopy colimit, and we denote
	\[\hoco(X_{\G})=\gmc(X_{\G}).\]
	
	By a \emph{map of diagrams} $f\colon X_{\G}\to Y_{\G}$ over the same graph $\G$, we mean a collection of maps $f_v\colon X_v\to Y_v$ for each vertex $v$. If every square of the form
	\[\begin{tikzcd}
		X_v\arrow{r}{\varphi_e}\arrow[swap]{d}{f_v} & X_w \arrow{d}{f_w}\\
		Y_v\arrow{r}{\psi_e} & Y_w
	\end{tikzcd}\]
	\begin{itemize}
		\item [a)] commutes, that is, $f_w\varphi_e=\psi_ef_v$, then we call $f$ a \emph{commutative} map of diagrams.
		\item [b)] commutes up to homotopy, that is, $f_w\varphi_e$ and $\psi_ef_v$ are homotopic, then we call $f$ a \emph{homotopy commutative} map of diagrams.
	\end{itemize}
	
	If a homotopy commutative map $f\colon X_{\G}\to Y_{\G}$ is given by homotopy equivalences, we refer to $X_{\G}$ and $Y_{\G}$ as \emph{homotopy equivalent diagrams}.
	
	The aforementioned fact that we require is a simple corollary of known theorems, and we defer its proof to the appendix.
	\begin{Lem}\label{simhoco}
		Homotopy equivalent source--drain diagrams have homotopy equivalent homotopy colimits.
	\end{Lem}
	
	\begin{Rem}\label{remhediag}
		Throughout the paper, we consider many diagrams built in the following way. Let $X_{\G}=(X_v,\varphi_e)_{\G}$ be a source--drain diagram of spaces. For each vertex $v$ of $\G$, consider a space $Y_v$ homotopy equivalent to $X_v$ via a homotopy equivalence $f_v\colon X_v\to Y_v$ with homotopy inverse $f'_v\colon Y_v\to X_v$.
		Let $Y_{\G}=(Y_v,\psi_e)_{\G}$ be a diagram given by spaces $Y_v$ and maps $\psi_e=f_w\circ\varphi_e\circ f'_v$. Then, $X_{\G}$ and $Y_{\G}$ are homotopy equivalent diagrams.
		
		In such instances we do not specify the maps in the diagrams. However, we have to assume that the choice of homotopy equivalences has been fixed. Different homotopy equivalences between $X_v$ and $Y_v$ yield the same homotopy type of homotopy colimits of $X_{\G}$ and $Y_{\G}$, but may induce non-homotopic homotopy equivalences between $\hoco(X_{\G})$ and $\hoco(Y_{\G})$. This is the reason why homotopy colimit is not the colimit in the homotopy category.
	\end{Rem}
	
	A frequently encountered example of a source--drain diagram is a pushout diagram:
	\[\begin{tikzcd}
		X & A \arrow[l, "f"']\arrow[r, "g"] & Y.
	\end{tikzcd}\]
	We denote the colimit of such a diagram by $X\sqcup_A Y$ and we call it the \emph{pushout} of this diagram. Similarly, we call the homotopy colimit of a pushout diagram the \emph{homotopy pushout}. The homotopy pushout is known as the double mapping cylinder.
	
	If $P$ is the homotopy pushout of the diagram above, then we also call the following diagram
	
	\[\begin{tikzcd}
		A \arrow{r}{g} \arrow[swap]{d}{f}\arrow{rd} & Y \arrow{d} \\
		X \arrow{r} & P  
	\end{tikzcd}\]
	a \emph{homotopy pushout diagram}. By Lemma~\ref{simhoco} and our definition of a homotopy colimit of a source--drain diagram, we also refer to any homotopy equivalent diagram as a homotopy pushout diagram.
	
	Homotopy pushout diagrams are useful as they are well-suited to the use of the Mayer--Vietoris sequence.
	
	\begin{Rem}\label{remnerve}
		We will refer to the proof of the Projection lemma (see \cite[Section~15.4]{kozlov2008combinatorial} or \cite[page~459]{hatcher2005algebraic}). This lemma states that for nerve diagrams associated with a covering of a space, a natural projection map $\hoco(X_{\G})\to\colim(X_{\G})$ is a homotopy equivalence. It concerns nerve diagrams, but its proof works as well in the cases we encounter, where the diagram comes from a covering, but is simplified by dropping redundant intersections.
	\end{Rem}
	
	\section{Deleted square of a cone}
	
	In \cite{wiltshiregordon2018configurationspaceproduct} Wiltshire-Gordon gives a description of the configuration space on the product of spaces in terms of generalized notion of the configuration spaces on factors, where the particles are allowed to collide in a specified way.
	
	The case of the deleted square (two particles) is easy to recall. $(X\times Y)^2_{\circ}$ can be decomposed into two subspaces:
	\[A=\{((x_1,y_1),(x_2,y_2))\in(X\times Y)^2_{\circ}\mid y_1\neq y_2 \}\simeq X^2\times Y^2_{\circ},\]
	and
	\[B=\{((x_1,y_1),(x_2,y_2))\in(X\times Y)^2_{\circ}\mid x_1\neq x_2 \}\simeq X^2_{\circ}\times Y^2,\]
	with intersection:
	\[A\cap B\simeq X^2_{\circ}\times Y^2_{\circ}.\]
	This yields a homotopy pushout diagram:
	\begin{equation}\label{gordon}
		\begin{tikzcd}
			X^2_{\circ}\times Y^2_{\circ} \arrow{r} \arrow{d} \arrow{dr}& X^2\times Y^2_{\circ} \arrow{d} \\
			X^2_{\circ}\times Y^2 \arrow{r} & (X\times Y)^2_{\circ}
		\end{tikzcd} 
	\end{equation}
	and applying the Mayer--Vietoris to it allows the description of homology of $(X\times Y)^2_{\circ}$ in simpler terms.
	
	We mimic this approach in case of $(\Cn X)^2_{\circ}$. We use the following decomposition:
	\[A=\{([x_1,t_1],[x_2,t_2])\in(\Cn X)^2_{\circ}\mid t_1\leq t_2, t_1 < 1\},\]
	and
	\[B=\{([x_1,t_1],[x_2,t_2])\in(\Cn X)^2_{\circ}\mid t_2\leq t_1, t_2 < 1\},\]
	to get a homotopy pushout diagram
	\begin{equation}\label{diag2}
		\begin{tikzcd}
			A\cap B \arrow{r} \arrow{d} \arrow{dr} & A \arrow{d} \\
			B \arrow{r} & (\Cn X)^2_{\circ}
		\end{tikzcd} 
	\end{equation}
	$(\Cn X)^2_{\circ}$ is not a CW-complex, and therefore we need to show that it is indeed a homotopy pushout of 
	\[\begin{tikzcd}
		A & A\cap B \arrow{l}\arrow{r} & B
	\end{tikzcd}\]
	by finding appropriate open neighbourhoods of $A$ and $B$ to use the Projection lemma. 
	
	We define open sets $U,V\subset (\Cn X)^2_{\circ}$ as follows:
	\[U=A\cup\left\{([x_1,t_1],[x_2,t_2])\in(\Cn X)^2_{\circ}\mid t_1< \min \{t_2 +\textstyle\frac{1}{2},1\}, x_1\neq x_2\right\},\]
	and
	\[V=B\cup\left\{([x_1,t_1],[x_2,t_2])\in(\Cn X)^2_{\circ}\mid t_2<\min \{t_1+\textstyle\frac{1}{2},1\}, x_1\neq x_2 \right\}.\]
	
	They are indeed open. For instance, the complement $U^c$ is a union of 
	two closed subsets given by $t_1\geq t_2+\frac{1}{2}$, and by $t_1=1$ and the third one given as an intersection of $(\Cn X)^2_{\circ}$ with a closed subset of $(\Cn X)^2$ given by $t_1\geq t_2\wedge x_1=x_2$.
	
	Moreover, $A$ is a deformation retract of $U$. If $t_1>t_2$, then $x_1\neq x_2$ and the particles $[x_1,t_1], [x_2,t_2]$ can be moved towards each other with respect to the second coordinate, leaving the first coordinate unchanged. In result the point $([x_1,t_1], [x_2,t_2])$ is moved towards $([x_1,\frac{t_1+t_2}{2}], [x_2,\frac{t_1+t_2}{2}])$ along the line segment connecting them in $(\Cn X)^2_{\circ}$. $V$ has an analogous deformation retraction onto $B$.
	
	The intersection $U\cap V$ is given by
	\[(|t_1-t_2|<\textstyle\frac{1}{2})\wedge (t_1,t_2<1)\wedge (x_1\neq x_2)\] and has a deformation retraction onto $A\cap B\cong X^2_{\circ}\times[0,1)$ which is constructed by moving the particles from the pair $([x_1,t_1],[x_2,t_2])$ to the subspace where the second coordinate is $\frac{t_1+t_2}{2}$ without changing the first coordinate --- it is a composition of deformation retractions of $U$ and $V$ described above.
	
	Projection lemma applied to the cover $\{U,V\}$ shows that $(\Cn X)^2_{\circ}$ is a homotopy pushout of 
	\[\begin{tikzcd}
		U & U\cap V \arrow{l}\arrow{r} & V,
	\end{tikzcd}\]
	and indicated deformation retractions allow the use of Remark~\ref{remhediag} and Lemma~\ref{simhoco} to prove the claim. $\square$
	
	Notice the homotopy equivalences:
	\[A\simeq X\text{ and }B\simeq X,\] which are obtained by pushing the particles apart with respect to the second coordinate so that one ends up in the base and the other in the apex of the cone; and
	\[A\cap B\simeq X^2_{\circ},\] which is obtained by pushing both particles to the base of the cone ($t_i=0$ subspace). By the construction from Remark~\ref{remhediag} and by Lemma~\ref{simhoco}, the following diagram is homotopy equivalent to the diagram in (\ref{diag2}).
	\begin{equation}\label{eq3}
		\begin{tikzcd}
			X^2_{\circ}\arrow[r, "p_1"] \arrow[d, "p_2"'] \arrow{dr} & X \arrow{d} \\
			X \arrow{r} & (\Cn X)^2_{\circ}
		\end{tikzcd} 
	\end{equation}
	where $p_i$ is a projection onto the $i$-th coordinate. That proves the lemma:
	\begin{Lem}\label{conepushout}
		$(\Cn X)^2_{\circ}$ is homotopy equivalent to a double mapping cylinder of projections of $X^2_{\circ}$ onto the coordinates.
	\end{Lem}
	\begin{Cor}\label{corA}
		If $X$ and $Y$ are both contractible, then \begin{enumerate}[label=\roman*)]
			\item $(X\times Y)^2_{\circ}\simeq X^2_{\circ}*Y^2_{\circ}$, and
			\item $(\Cn X)^2_{\circ}\simeq S(X^2_{\circ})$.
		\end{enumerate}
	\end{Cor}
	\begin{proof}
		Diagram~(\ref{gordon}) is in this case homotopy equivalent to 
		
		\[\begin{tikzcd}
			X^2_{\circ}\times Y^2_{\circ} \arrow{r} \arrow{d} \arrow{dr} & Y^2_{\circ} \arrow{d} \\
			X^2_{\circ} \arrow{r} & X^2_{\circ}*Y^2_{\circ}
		\end{tikzcd}\]
		and diagram~(\ref{eq3}) is homotopy equivalent to
		\[\begin{tikzcd}
			X^2_{\circ}\arrow{r}{p_1} \arrow{d}{p_2} \arrow{dr} & * \arrow{d} \\
			* \arrow{r} & S(X^2_{\circ})
		\end{tikzcd}\]
		
	\end{proof}
	This gives a non-obvious example of non-homeomorphic spaces with homotopy equivalent deleted squares:
	\begin{Cor}
		For contractible $X$: $(X\times I)^2_{\circ}\simeq(\Cn X)^2_{\circ}.$
	\end{Cor}

	\subsection{Example: deleted square of a cone over a star graph}
	A star graph $St(n)$ is a tree with $n+1$ vertices from which $n$ are leaves (are of valence 1). We examine the case $n=3$.
	
	$(\Cn(St_3))^2_{\triangle}$ is homeomorphic to a sphere. It can be presented as a subdivision of the surface of a cube. Figure~\ref{fig:Y} shows the 1-skeleton of $(\Cn(St_3))^2_{\triangle}$. Its vertices are pairs of vertices of $St_3$, which are identified with numbers $0,1,2,3$, where nonzero ones correspond to leaves, and $a$ is the apex of the cone.
	
	There are two subcomplexes of $(\Cn(St_3))^2_{\triangle}$ isomorphic to $St_3$ --- they correspond to one of the particles occupying the apex; and there is a subcomplex isomorphic to $(St_3)^2_{\triangle}$ that correspond to both particles lying in the base of the cone. All of them are indicated by thick lines to visualize Lemma~\ref{conepushout}.
	\begin{figure}
		\centering
		{\begin{tikzpicture}[x={(1.2cm,-0.3cm)}, 
				y={(0.8cm,0.5cm)}, 
				z={(0cm,1.3cm)}, scale=1, line join=round, line cap=round]
				\coordinate (A) at (0,0,0);
				\coordinate (B) at (3,0,0);
				\coordinate (C) at (3,3,0);
				\coordinate (D) at (0,3,0);
				
				\coordinate (E) at (0,0,3);
				\coordinate (F) at (3,0,3);
				\coordinate (G) at (3,3,3);
				\coordinate (H) at (0,3,3);
				
				\coordinate (AB) at ($(A)!0.5!(B)$);
				\coordinate (BC) at ($(B)!0.5!(C)$);
				\coordinate (CD) at ($(C)!0.5!(D)$);
				\coordinate (AD) at ($(A)!0.5!(D)$);
				
				\coordinate (EF) at ($(E)!0.5!(F)$);
				\coordinate (FG) at ($(F)!0.5!(G)$);
				\coordinate (GH) at ($(G)!0.5!(H)$);
				\coordinate (EH) at ($(E)!0.5!(H)$);
				
				\coordinate (AE) at ($(A)!0.5!(E)$);
				\coordinate (BF) at ($(B)!0.5!(F)$);
				\coordinate (CG) at ($(C)!0.5!(G)$);
				\coordinate (DH) at ($(D)!0.5!(H)$);
				
				\draw (A) -- (B) -- (C) -- (D) -- cycle;
				\draw (E) -- (F) -- (G) -- (H) -- cycle;
				\draw (A) -- (E) (B) -- (F) (C) -- (G) (D) -- (H);
				
				\draw (H) -- (A) (H) -- (F) (H) -- (C);
				\draw (B) -- (D) (B) -- (E) (B) -- (G);
				\draw (AB) -- (AD) (BC) -- (CD) (BC) -- (CG) (BF) -- (FG);
				\draw (AB) -- (AE) (BF) -- (EF) (AE) -- (EH) (AD) -- (DH);
				\draw (CD) -- (DH) (CG) -- (GH) (EF) -- (EH) (FG) -- (GH);
				
				\draw[very thick] (AB) -- (B) -- (BF);
				\draw[very thick] (B) -- (BC);
				\draw[very thick] (DH) -- (H) -- (GH);
				\draw[very thick] (H) -- (EH);
				
				\draw[very thick] (A) -- (AD) -- (D) -- (CD) -- (C) -- (CG) -- (G) -- (FG) -- (F) -- (EF) -- (E) -- (AE) -- cycle;
				
				\foreach \v in {A,C,D,E,F,G,AD,CD,CG,FG,EF,AE,B,H,AB,BC,BF,DH,EH,GH} {
					\fill[black] (\v) circle [radius=2pt];
				}
				\tikzset{label style/.style={font=\scriptsize, inner sep=2pt, fill=white, fill opacity=0.6, text opacity=1}}
				
				\node[label style, below left]  at (A) {$(1,0)$};
				\node[label style, below right] at (B) {$(0,a)$};
				\node[label style, below right]       at (C) {$(2,0)$};
				\node[label style, above left]  at (D) {$(0,3)$};
				\node[label style, left]        at (E) {$(0,2)$};
				\node[label style, below right] at (F) {$(3,0)$};
				\node[label style, above right] at (G) {$(0,1)$};
				\node[label style, above]       at (H) {$(a,0)$};
				
				\node[label style, below left]       at (AB) {$(1,a)$};
				\node[label style, below right]       at (BC) {$(2,a)$};
				\node[label style, above right] at (CD) {$(2,3)$};
				\node[label style, below right]        at (AD) {$(1,3)$};
				
				\node[label style, below left]       at (EF) {$(3,2)$};
				\node[label style, below right]       at (FG) {$(3,1)$};
				\node[label style, above]       at (GH) {$(a,1)$};
				\node[label style, left]        at (EH) {$(a,2)$};
				
				\node[label style, left]        at (AE) {$(1,2)$};
				\node[label style, below left]       at (BF) {$(3,a)$};
				\node[label style, right]       at (CG) {$(2,1)$};
				\node[label style, above right]        at (DH) {$(a,3)$};
		\end{tikzpicture}}
		\caption{1-skeleton of $(\Cn(St_3))^2_{\triangle}$}
		\label{fig:Y}
	\end{figure}
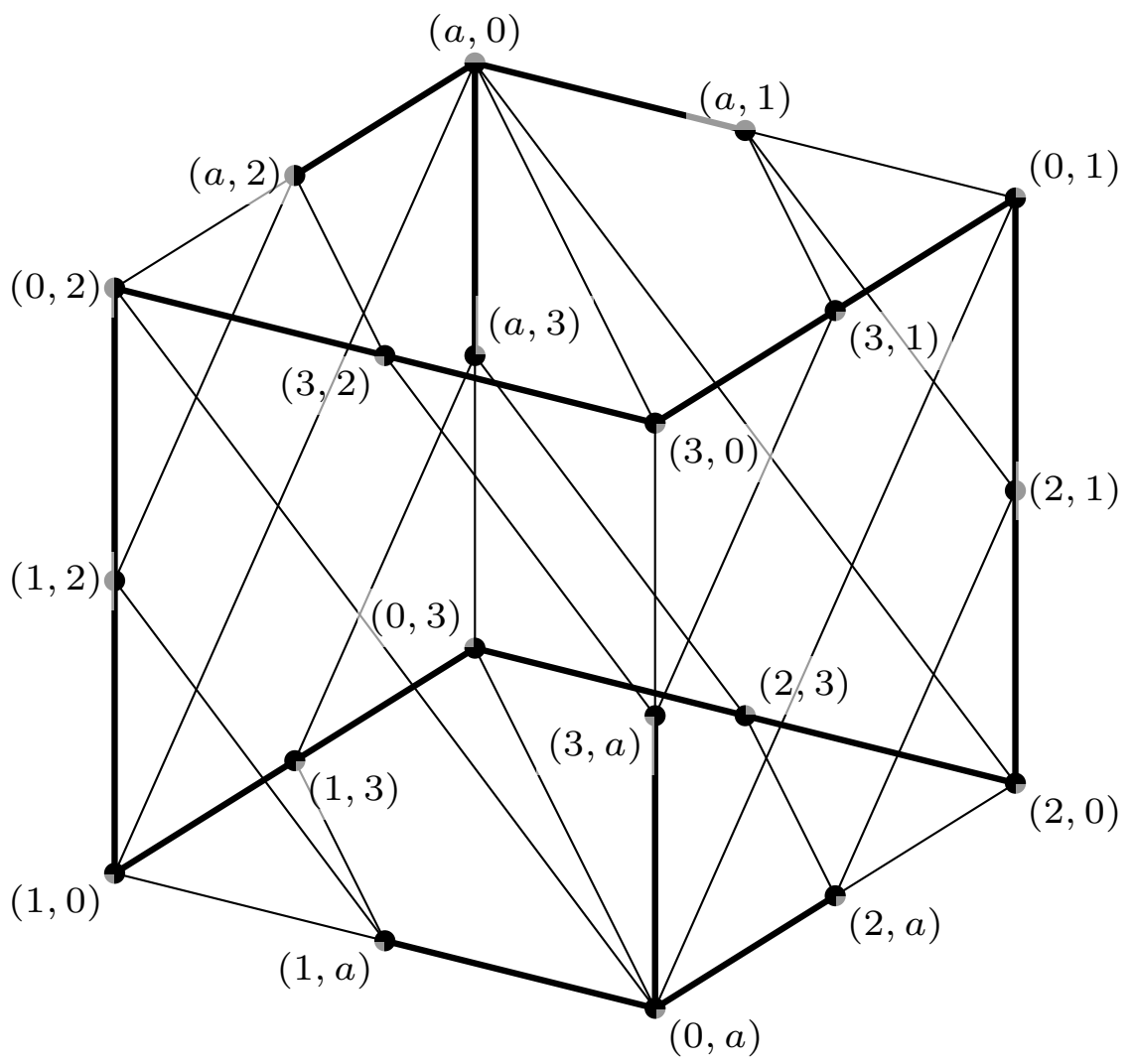
	\begin{Rem}
		One can also see $St(4)^2_{\triangle}$ in Figure~\ref{fig:Y}. Indeed, $St(4)$ can be identified with subcomplex of $(\Cn(St_3)$ consisting of the base and the edge connecting the vertex $0$ to the apex $a$ of the cone. $St(4)^2_{\triangle}$ is a subdivision of the $1$-skeleton of a cube.
	\end{Rem}
	\subsection{Simplicial deleted join}
	The deleted squares of cones have been investigated in \cite{copeland1965deleted} where cohomology of $(\Cn X)^2_{\triangle}$ is described and $(\Cn X)^2_{\triangle}$ is shown to be homeomorphic with the subspace of the join $X*X$ given as the image of \[(X^2_{\triangle}\times [-1,1])\cup (X\times\{-1,1\})\] through the quotient map. This subspace is called the \emph{simplicial deleted join} of $X$ with itself.
	
	We define simplicial deleted join combinatorially as:
	\[X\sdj Y=\{(\sigma\times\{-1\})\cup(\tau\times\{1\})\mid\sigma\in X,\tau\in Y,\sigma\cap\tau=\emptyset\}.\]
	Its geometric realization is homeomorphic with the subspace of $X*Y$ given by \[((X\sdp Y)\times (-1,1))\cup (X\times\{-1\})\cup (Y\times\{1\}),\] but we omit the proof. We denote $X^{2*}_{\triangle}:=X\sdj X$.
	
	To see that $(\Cn X)^2_{\triangle}$ is homeomorphic to $X^{2*}_{\triangle}$ consider a map $(\Cn X)^2_{\triangle}\to X^{2*}_{\triangle}$ given by
	\[([x,t],[y,s])\mapsto[x,y,t-s].\]
	
\begin{figure}
		\centering
		\begin{subfigure}{.5\textwidth}
			\centering
			\begin{tikzpicture}[scale=0.6, line join=round, line cap=round]
				\coordinate (A) at (-0.66,0.39);
				\coordinate (B) at (0.66,0.39);
				\coordinate (C) at (0,-0.77);
				
				\coordinate (A') at (3,-1.74);
				\coordinate (B') at (-3,-1.74);
				\coordinate (C') at (0,3.48);
				
				\draw (A) -- (B) -- (C) -- cycle;
				\draw (A') -- (B') -- (C') -- cycle;
				\draw (A') -- (B) -- (C') -- (A) -- (B') -- (C) -- cycle;
				
				\fill[opacity=0.2, even odd rule] (A) -- (B) -- (C) -- cycle (A') -- (B') -- (C') -- cycle;
				
				\coordinate (AB') at ($(A)!0.5!(B')$);
				\coordinate (B'C) at ($(B')!0.5!(C)$);
				\coordinate (CA') at ($(C)!0.5!(A')$);
				\coordinate (A'B) at ($(A')!0.5!(B)$);
				\coordinate (BC') at ($(B)!0.5!(C')$);
				\coordinate (C'A) at ($(C')!0.5!(A)$);
				
				\draw (A'B) -- (BC') -- (C'A) -- (AB') -- (B'C) -- (CA') -- cycle;
			\end{tikzpicture}
			\caption{$(\Cn(\partial\Delta^2))^{2}_{\triangle}$}
		\end{subfigure}%
		\begin{subfigure}{.5\textwidth}
			\centering
			\begin{tikzpicture}[scale=0.6, line join=round, line cap=round]
				\coordinate (A) at (-0.66,0.39);
				\coordinate (B) at (0.66,0.39);
				\coordinate (C) at (0,-0.77);
				
				\coordinate (A') at (3,-1.74);
				\coordinate (B') at (-3,-1.74);
				\coordinate (C') at (0,3.48);
				
				\draw (A) -- (B) -- (C) -- cycle;
				\draw (A') -- (B') -- (C') -- cycle;
				\draw (A') -- (B) -- (C') -- (A) -- (B') -- (C) -- cycle;
				
				\fill[opacity=0.2, even odd rule] (A) -- (B) -- (C) -- cycle (A') -- (B') -- (C') -- cycle;	
			\end{tikzpicture}
			\caption{$(\partial\Delta^2)^{2*}_{\triangle}$}
		\end{subfigure}
		\caption{Comparison between simplicial deleted square of a cone over $\partial\Delta^2$ and a simplicial deleted join of $\partial\Delta^2$ with itself.}
		\label{fig:tri}
	\end{figure}
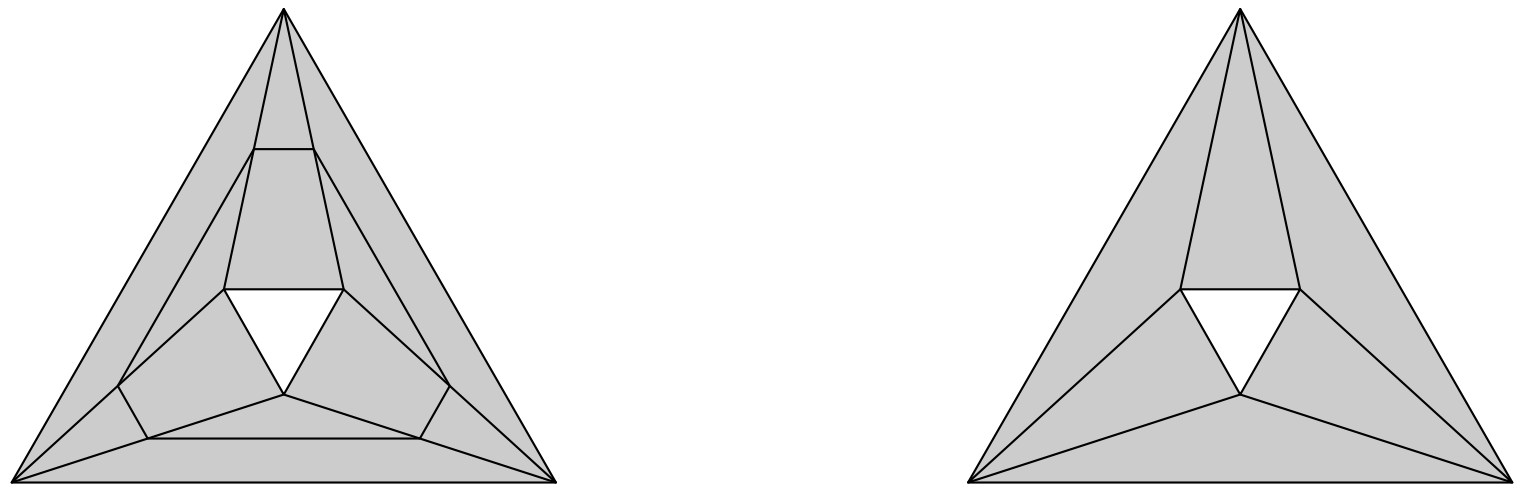

	In particular, $(\Cn X)^2_{\triangle}$ is a subdivision of $X^{2*}_{\triangle}$ (see Figure~\ref{fig:tri}).
	
	While it is relevant for the proofs, we will write $(\Cn X)^2_{\triangle}$, but later we will use $X^{2*}_{\triangle}$ to simplify notation.
	
	\subsection{Cone over a graph}
	Let $X=\G$ be a connected graph. Apply the Mayer--Vietoris sequence to the diagram~(\ref{eq3}) to get:
	
	\begin{equation}\label{LEScone}
		\begin{tikzcd}
			0 \arrow{r} & H_3((\Cn\G)^2_{\circ})\arrow{r} & H_2(\G^2_{\circ})\arrow{r} & 0 \arrow{r} & \text{\color{white}{a}} \\
			\arrow{r} & H_2((\Cn\G)^2_{\circ})\arrow{r} & H_1(\G^2_{\circ})\arrow{r} &  H_1(\G)\oplus H_1(\G)\arrow{r} & \text{\color{white}{a}}\\
			\arrow{r} & H_1((\Cn\G)^2_{\circ})\arrow{r} & \Z\arrow{r} &\Z\oplus\Z\arrow{r} &\text{\color{white}{a}}\\
			\arrow{r} & \Z\arrow{r} & 0
		\end{tikzcd} 
	\end{equation}
	Proposition 1.3 in \cite{barnett2009topology} states that the map
	$H_1(\G^2_{\circ})\to H_1(\G)\oplus H_1(\G)$
	is an epimorphism if $\G$ is not homeomorphic to a circle (it is also implied by Theorem~\ref{simconprop}). This completes the proof of
	\begin{Prop}[{Copeland \cite{copeland1965deleted}}]
		Let $\G$ be a nondegenerate graph. Then
		\begin{enumerate}[label=\roman*)]
			\item $H_2(\G^2_{\circ})\cong H_3((\Cn\G)^2_{\circ})$.
			\item $H_1(\G^2_{\circ})\cong H_2((\Cn\G)^2_{\circ})\oplus H_1(\G)\oplus H_1(\G)$.
		\end{enumerate}
	\end{Prop}
	In particular, the Betti numbers of $\G^2_{\circ}$ can be computed knowing the Betti numbers of $(\Cn\G)^2_{\circ}$ and vice versa.
	\subsection{Example: complete bipartite graphs}
	
	Let $K_{n,m}$ denote a complete bipartite graph. In \cite{farber2010topology}, the authors computed the homology of $(K_{n,m})^2_{\circ}$. We compute the homotopy type of $(\Cn K_{n,m})^2_{\circ}$.
	
	If $\overline n$ denotes a set of $n$ points with discrete topology, then $K_{n,m}\cong\overline n*\overline m$. Thus,
	\[\Cn K_{n,m}\cong\Cn(\overline n*\overline m)\cong\Cn(\overline n)\times\Cn(\overline m),\]
	and $\Cn(\overline n)$ is just a star graph $St_n$.
	
	By Corollary~\ref{corA}, $(St_n\times St_m)^2_{\circ}\simeq(St_n)^2_{\circ}*(St_m)^2_{\circ}$. Thus,
	\[(\Cn K_{n,m})^2_{\circ}\simeq(St_n)^2_{\circ}*(St_m)^2_{\circ}.\] 
	Notice that $St_1\cong St_2\cong[0,1]$, therefore, $(St_1)^2_{\circ}\simeq(St_2)^2_{\circ}\simeq\overline 2$, and
	\[(\Cn K_{1,m})^2_{\circ}\simeq(\Cn K_{2,m})^2_{\circ}\simeq S((St_m)^2_{\circ}).\]
	For $n\geq 3$ we know that $(St_n)^2_{\circ}$ has a homotopy type of bouquet of $n^2-3n+1$ circles (see Lemma~\ref{trees}). Thus, $(\Cn K_{n,m})^2_{\circ}$ has the homotopy type of a bouquet of $(n^2-3n+1)(m^2-3m+1)$ spheres $S^3$, and we get
	\begin{Prop}
		Complete bipartite graphs $K_{n,m}$
		\begin{enumerate}[label=\alph*)]
			\item satisfy $(\Cn K_{n,m})^2_{\circ}\simeq S^1$ for $n,m<3$,
			\item are unripe for $\min\{n,m\}<3$ and $\max\{n,m\}\geq3$, and
			\item are ripe for $n,m\geq 3$.
		\end{enumerate}
		Moreover, for 
		unripe 
		$K_{n,m}$:
		\[(\Cn K_{n,m})^2_{\circ}\simeq\bigvee_{i=1}^{k^2-3k+1}S^2,\]
		where $k=\max\{n,m\}$, and for ripe $K_{n,m}$:
		\[(\Cn K_{n,m})^2_{\circ}\simeq\bigvee_{i=1}^{(n^2-3n+1)(m^2-3m+1)}S^3.\]
	\end{Prop}
	\section{Main theorem}
	Recall the graph enlarging operations from the introduction:
	\begin{itemize}
		\item \emph{Adding a new vertex as a leaf}: for a graph $\G$ and its vertex $v$, we define a new graph $\adv(\G,v)$, which is obtained by identifying one endpoint of $[0,1]$ with~$v$.
		\item \emph{Adding a new edge between existing vertices}: for a graph $\G$ and its distinct vertices $v_1$ and $v_2$, we define a new graph $\ade(\G,v_1,v_2)$, which is obtained by identifying the endpoints of $[0,1]$ with $v_1$ and $v_2$. 
	\end{itemize}
	
	\begin{Thm}\label{main}
		For a graph $\G$ and its vertices $v$, $v_1$, and $v_2$, with $v_1\neq v_2$:
		\begin{enumerate}[label=\alph*)]
			\item $(\Cn(\adv(\G,v)))^2_{\triangle}$ is a homotopy colimit of the following diagram:
			\[\begin{tikzcd}[column sep=1.5 em]
				* & \Cn\G\sdp\Cn\{v\} \arrow{l}\arrow{r} & (\Cn\G)^2_{\triangle} & \Cn \{v\}\sdp\Cn\G \arrow{l}\arrow{r} & *,
			\end{tikzcd}\]
			where
			\[\Cn\G\sdp\Cn\{v\}\simeq \bigvee_{i=1}^{\nu(v)-1}S^1\simeq \Cn \{v\}\sdp\Cn\G.\]
			\item $(\Cn(\ade(\G,v_1,v_2)))^2_{\triangle}$ is a homotopy colimit of the following diagram:
			\[\begin{tikzcd}[column sep=1.3 em]
				* & \Cn\G\sdp\Cn\{v_1,v_2\} \arrow{l}\arrow{r} & (\Cn\G)^2_{\triangle} & \Cn \{v_1,v_2\}\sdp\Cn\G \arrow{l}\arrow{r} & *,
			\end{tikzcd}\]
			where \[\Cn\G\sdp\Cn\{v_1,v_2\}\simeq \left(\bigvee_{i=1}^{b_0-1}S^1\right)\vee\left(\bigvee_{i=1}^{b_1}S^2\right)\simeq \Cn \{v_1,v_2\}\sdp\Cn\G,\]
			and $b_0$ and $b_1$ are the zeroth and the first Betti numbers of $\G\setminus\{v_1,v_2\}$.
		\end{enumerate}
	\end{Thm}
	By Lemma~\ref{simhoco}, we can formulate this theorem more concisely:
	\begin{Cor}\label{maincor}
		For a graph $\G$ and its vertices $v$, $v_1$, and $v_2$, with $v_1\neq v_2$:
		\begin{enumerate}[label=\alph*)]
			\item\label{maincora} $(\Cn(\adv(\G,v)))^2_{\triangle}$ is a homotopy colimit of the following diagram:
			\[\begin{tikzcd}[column sep=1.5 em]
				* & \bigvee_{i=1}^{\nu(v)-1}S^1 \arrow{l}\arrow{r} & (\Cn\G)^2_{\triangle} & \bigvee_{i=1}^{\nu(v)-1}S^1 \arrow{l}\arrow{r} & *.
			\end{tikzcd}\]
			\item\label{maincorb} $(\Cn(\ade(\G,v_1,v_2)))^2_{\triangle}$ is a homotopy colimit of the following diagram:
			\[\begin{tikzcd}[column sep=1.2 em]
				* & \left(\bigvee_{i=1}^{b_0-1}S^1\right)\vee\left(\bigvee_{i=1}^{b_1}S^2\right) \arrow{l}\arrow{r} & (\Cn\G)^2_{\triangle} & \left(\bigvee_{i=1}^{b_0-1}S^1\right)\vee\left(\bigvee_{i=1}^{b_1}S^2\right) \arrow{l}\arrow{r} & *,
			\end{tikzcd}\]
		\end{enumerate}
		where $b_0$ and $b_1$ are the zeroth and the first Betti numbers of $\G\setminus\{v_1,v_2\}$.
	\end{Cor}
	According to Lemma~\ref{lemret}, we restrict our attention to the simplicial deleted square, bearing in mind that it has the same homotopy type as the usual deleted square.
	
	We prepare our tools in a general manner. Suppose that a simplicial complex $X$ is a union of subcomplexes $A$ and $B$ with intersection $A\cap B= C$.
	Then, there is a decomposition:
	\[X^2_{\triangle}=A^2_{\triangle}\cup B^2_{\triangle}\cup (A\sdp B)\cup (B\sdp A).\]
	Notice also that the simplicial deleted product is distributive with respect to the intersection:
	\[(X_1\sdp Y_1)\cap (X_2 \sdp Y_2)=(X_1\cap X_2)\sdp(Y_1\cap Y_2).\]
	Therefore:
	\begin{itemize}
		\item $A^2_{\triangle}\cap (A\sdp B) = A\sdp C$ and $A^2_{\triangle}\cap (B\sdp A) = C\sdp A$,
		\item $B^2_{\triangle}\cap (A\sdp B) = C\sdp B$ and $B^2_{\triangle}\cap (B\sdp A) = B\sdp C$,
		\item $A^2_{\triangle}\cap B^2_{\triangle}\cap (A\sdp B) \cap (B\sdp A) = C^2_{\triangle}$.
	\end{itemize}
	Consider the following diagram:
	
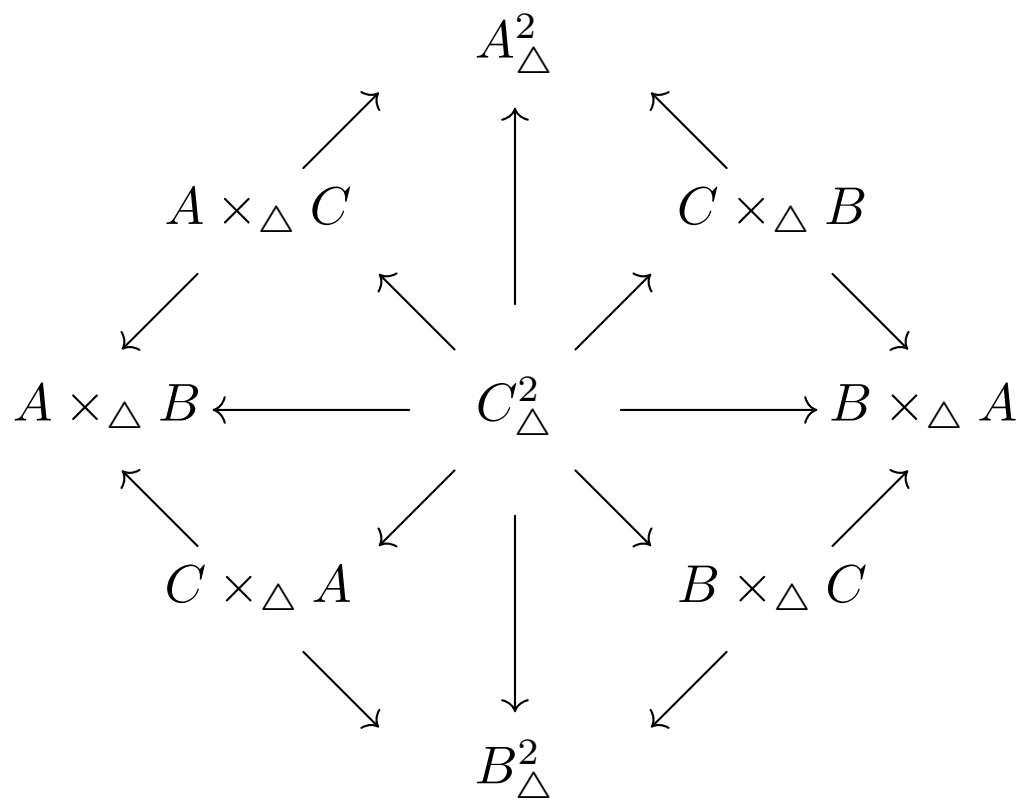
\begin{center}
		\begin{tikzpicture}[line cap=round,line join=round,x=1.0cm,y=1.0cm]
			\clip (-3.5,-3) rectangle (3.5,3);
			
			\draw [->] (0.7,0)--(2,0);
			\draw [->] (-0.7,0)--(-2,0);
			\draw [->] (0,0.7)--(0,2);
			\draw [->] (0,-0.7)--(0,-2);
			
			\draw [->] (-2.1,-0.9)--(-2.6,-0.4);
			\draw [->] (-1.4,-1.6)--(-0.9,-2.1);
			\draw [->] (2.1,-0.9)--(2.6,-0.4);
			\draw [->] (1.4,-1.6)--(0.9,-2.1);
			\draw [->] (-2.1,0.9)--(-2.6,0.4);
			\draw [->] (-1.4,1.6)--(-0.9,2.1);
			\draw [->] (2.1,0.9)--(2.6,0.4);
			\draw [->] (1.4,1.6)--(0.9,2.1);
			
			\draw [->] (-0.4,-0.4)--(-0.9,-0.9);
			\draw [->] (-0.4,0.4)--(-0.9,0.9);
			\draw [->] (0.4,-0.4)--(0.9,-0.9);
			\draw [->] (0.4,0.4)--(0.9,0.9);
			
			\draw (-2.7,0) node {$A \sdp B$};
			\draw (2.7,0) node {$B \sdp A$};
			\draw (0,2.4) node {$A^2_{\triangle}$};
			\draw (0,-2.4) node {$B^2_{\triangle}$};
			
			\draw (-1.7,1.3) node {$A \sdp C$};
			\draw (1.7,1.3) node {$C \sdp B$};
			\draw (-1.7,-1.2) node {$C \sdp A$};
			\draw (1.7,-1.2) node {$B \sdp C$};
			\draw (0,0) node {$C^2_{\triangle}$};
		\end{tikzpicture} 
	\end{center}
	where all maps are inclusions.
	
	It is not a nerve diagram, but $X^2_{\triangle}$ is homeomorphic to its colimit, and the proof of the Projection lemma applies to it (Remark~\ref{remnerve}). Thus, if homotopy types of indicated parts are known, we can obtain the homotopy type of $X^2_{\triangle}$ by using the general homotopy colimit. However, a source--drain diagram can be obtained in the following particularly useful case.
	
	Let $B\subset A$ be a subcomplex and consider $(A\sqcup_B \Cn B)^2_{\triangle}$. From the definition of the simplicial deleted square, we get:
	\[(\Cn B)^2_{\triangle}\subset \left(A\sdp\Cn B\right)\cup \left(\Cn B\sdp A\right).\]
	Hence, 
	\[(A\sqcup_B \Cn B)^2_{\triangle}=A^2_{\triangle}\cup\left(A\sdp\Cn B\right)\cup \left(\Cn B\sdp A\right).\]
	Moreover, 
	\[\left(A\sdp\Cn(B)\right)\cap \left(\Cn(B)\sdp A\right)=A\sdp B,\]
	and $A\sdp B\subset A^2_{\triangle}$. Thus, $(A\sqcup_B \Cn B)^2_{\triangle}$ is homeomorphic to the colimit of the diagram
	\begin{equation}\label{diag0}
		\begin{tikzcd}
			A\sdp\Cn B & A\sdp B \arrow{l}\arrow{r} & A^2_{\triangle} & B\sdp A\arrow{l}\arrow{r} & \Cn B\sdp A,
		\end{tikzcd}
	\end{equation}
	and by Remark~\ref{remnerve}, $(A\sqcup_B \Cn B)^2_{\triangle}$ is a homotopy colimit of this diagram.
	
	\begin{Rem}\label{remdcn}
		$A\sdp\Cn(B)$ is homotopy equivalent to $A$ as it deformation retracts to $A\times\{a\}$, where $a$ is the apex of a cone over $B$, by moving the second particle to $a$. It is important to keep in mind that $a$ is not a point of $A$ here, so that there will be no confusion later, when $A$ is itself a cone over some space $X$ and $B$ is a subspace of $A$ that is a cone over a subspace $Y\subset X$ (in such case we will refer to $B$ as \emph{subcone} of $X$). The same holds for inverted coordinates.
	\end{Rem}
	We implement this to the case of a cone. Let $Y\subset X$ be a subcomplex and notice that 
	\[\Cn(X\sqcup_Y\Cn Y)=\Cn X\sqcup_{\Cn Y}\Cn(\Cn Y),\]
	where $\Cn Y$ on the right is a subcone of $\Cn X$. Moreover, from Remark~\ref{remdcn},
	\[\Cn X \sdp \Cn(\Cn Y) \simeq * \simeq \Cn(\Cn Y) \sdp \Cn X.\]
	Thus, the diagram~(\ref{diag0}) is homotopy equivalent to:
	\begin{equation}\label{diag1}
		\begin{tikzcd}
			* & \Cn X \sdp \Cn Y \arrow{l}\arrow{r} & (\Cn X)^2_{\triangle} & \Cn Y\sdp\Cn X \arrow{l}\arrow{r} & *.
		\end{tikzcd}
	\end{equation}
	Notice that Remark~\ref{remdcn} does not apply to $\Cn X\sdp\Cn Y$ and $\Cn Y\sdp\Cn X$ here.
	
	We have proven the following lemma.
	
	\begin{Lem}\label{diaglem}
		$\Cn(X\sqcup_Y\Cn Y)^2_{\triangle}$ is a homotopy colimit of the diagram
		\[\begin{tikzcd}
			* & \Cn X \sdp \Cn Y \arrow{l}\arrow{r} & (\Cn X)^2_{\triangle} & \Cn Y\sdp\Cn X \arrow{l}\arrow{r} & *.
		\end{tikzcd}\]
	\end{Lem}
	
	The operations $\adv, \ade$ of adding a vertex and adding an edge can be realized at the level of topological spaces as:
	\[\adv(\G,v)=\G\sqcup_{\{v\}}\Cn\{v\},\]
	and
	\[\ade(\G,v_1, v_2)=\G\sqcup_{\{v_1,v_2\}}\Cn\{v_1,v_2\},\]
	where $\{v\}$ and $\{v_1,v_2\}$ are to be considered as discrete subspaces.
	
	Thus, by Lemma~\ref{diaglem}, $(\Cn(\adv(\G,v)))^2_{\triangle}$ and $(\Cn(\ade(\G,v_1,v_2)))^2_{\triangle}$ are homotopy colimits of the respective diagrams:
	\begin{equation}\label{diagL}
		\begin{tikzcd}[column sep=1.3 em]
			* & \Cn\G\sdp\Cn\{v\} \arrow{l}\arrow{r} & (\Cn\G)^2_{\triangle} & \Cn \{v\}\sdp\Cn\G \arrow{l}\arrow{r} & *
			\end{tikzcd}
	\end{equation}
	and
	\begin{equation}\label{diagE}
		\begin{tikzcd}[column sep=1.3 em]
			* & \Cn\G\sdp\Cn\{v_1,v_2\} \arrow{l}\arrow{r} & (\Cn\G)^2_{\triangle} & \Cn \{v_1,v_2\}\sdp\Cn\G \arrow{l}\arrow{r} & *
		\end{tikzcd}
	\end{equation}
	\begin{Lem}\label{pieces}
		Let $\G$ be a graph, and let $v$, $v_1$, and $v_2$ be its vertices, with $v_1\neq v_2$. Then,
		\begin{enumerate}[label=\alph*)]
			\item $\Cn\G\sdp\Cn\{v\}$ is homotopy equivalent to $\bigvee_{i=1}^{\nu(v)-1}S^1$.
			\item $\Cn\G\sdp\Cn\{v_1,v_2\}$ is homotopy equivalent to $S(\G\setminus\{v_1,v_2\})$.
		\end{enumerate}
		Analogous statements hold for inverted coordinates, i.e., for $\Cn\{v\}\sdp\Cn\G$ and $\Cn\{v_1,v_2\}\sdp\Cn\G$.
	\end{Lem}
	\begin{proof}
		We first analyze the more complex case in detail.
		
		Recall that the \emph{open star} $\st(v)$ of a vertex $v$ in a simplicial complex $X$ is the union of interiors of all simplices containing $v$ and $v$ itself, and the \emph{link} of $v$ is defined as the topological boundary of $st(v)$, that is, $\Lk(v)=\overline{st(v)}\setminus \st(v)$.
		
		Notice that $\Cn\G\sdp\Cn\{v_1,v_2\}$ is homeomorphic to the homotopy colimit of the diagram:
		\[\Cn(\G\setminus\st(v_1))\xleftarrow{}\G\setminus\st(v_1)\to\G\xleftarrow{}\G\setminus\st(v_2)\to\Cn(\G\setminus\st(v_2)),\]
		in which we display the homeomorphism classes according to the following homeomorphisms:
		\begin{itemize}
			\item $\Cn(\G\setminus\st(v_1))\cong\Cn(\G\setminus\st(v_1))\times\{v_1\}\cong\Cn\G\sdp\{v_1\}$
			\item $\G\setminus\st(v_1)\times(0,1)\cong\G\setminus\st(v_1)\times(v_1,a)\cong\Cn\G\sdp(v_1,a)$
			\item $\G\cong\G\times\{a\}\cong\Cn\G\sdp\{a\}$
		\end{itemize}
		and analogously for $v_2$ in place of $v_1$. Here by $(v_1, a)$ we mean the interior of the edge between $v_1$ and $a$. The first space corresponds to the second particle lying at $v_1$; the second to the second particle lying in the interior of the edge between vertex $v_1$ and the apex $a$ of $\Cn\{v_1,v_2\}$; and the third to the second particle lying at the apex of $\Cn\{v_1,v_2\}$. We abuse notation in the second bullet when using $\sdp$, but it should be clear from the context what we mean.
		
		Also notice that the resulting space is homotopy equivalent to the one obtained by coning $\G$ over $\G\setminus \st(v_1)$ and over $\G\setminus \st(v_2)$.
		
		If the open stars $\st(v_1)$, $\st(v_2)$ have empty intersection, then there is a deformation retraction of this space that maps $\Cn(\st(v_i))$ to $\Cn(\Lk(v_i))$ by pushing the bases into the cones, and preserves the rest of the space. The deformation retract is just the suspension $S(\G\setminus(\st(v_1)\cup\st(v_2))$.
		
		If $\st(v_1)$ and $\st(v_2)$ have nonempty intersection, it means, that there already exists an edge between $v_1$ and $v_2$, but we can subdivide all the edges between $v_1$ and $v_2$, without changing the homeomorphism class, making $\st(v_1)$ and $\st(v_2)$ disjoint.
		
		In both cases, the resulting space is homotopy equivalent to $S(\G\setminus\{v_1,v_2\})$, even if $\G$ here is not subdivided, and that is why we chose this to formulate the statement of the lemma.
		
		To prove (a) we proceed analogously. $\Cn\G\sdp\Cn\{v\}$ is given as the homotopy colimit of the diagram 
		\[\begin{tikzcd}
			\Cn(\G\setminus\st(v)) & \G\setminus\st(v)\arrow{l}\arrow{r} & \G.
		\end{tikzcd}\]
		Hence, it is homotopy equivalent to the mapping cone of the inclusion $\G\setminus\st(v)\rightarrow\G$. Thus, \[\Cn\G\sdp\Cn\{v\}\simeq \quot{\G}{\G\setminus\st(v)}\simeq \bigvee_{i=1}^{\nu(v)-1}S^1.\]
		By mapping cone of $f\colon X\to Y$ we mean the homotopy colimit of the diagram
		\[\begin{tikzcd}
			* & X\arrow{l}\arrow{r}{f} & Y.
		\end{tikzcd}\]
	\end{proof}
	To complete the proof of the main theorem just notice that 
	\[\Cn\G\sdp\Cn\{v\}\simeq\Cn \{v\}\sdp\Cn\G,\]
	and
	\[\Cn\G\sdp\Cn\{v_1,v_2\}\simeq\Cn \{v_1,v_2\}\sdp\Cn\G,\]
	and also
	\[S(\G\setminus\{v_1,v_2\})\simeq\left(\bigvee_{i=1}^{b_0-1}S^1\right)\vee\left(\bigvee_{i=1}^{b_1}S^2\right),\]
	where $b_i=\text{rk }H_i(\G\setminus\{v_1,v_2\})$.
	
	\subsection{Deleted join}
	Viewing the theorem from the deleted join perspective provides more geometric intuition. Notice that $(\ade(\G,v_1,v_2))^{2*}_{\triangle}$ is homeomorphic to the colimit of the diagram
	\[\begin{tikzcd}[column sep=1.5 em]
		\G\sdj [v_1,v_2] & \G\sdj \{v_1,v_2\}\arrow{l}\arrow{r} & \G^{2*}_{\triangle} & \{v_1,v_2\}\sdj\G \arrow{l}\arrow{r} & {[v_1,v_2]} \sdj \G.
	\end{tikzcd}\]
	Parts of this decomposition correspond to the spaces in the diagram~(\ref{diagE}):
	\begin{itemize}
		\item $*\simeq\Cn\G\sdp\Cn(\Cn\{v_1,v_2\})\cong\G\sdj [v_1,v_2]$,
		\item $\Cn\G\sdp\Cn\{v_1,v_2\}\cong\G\sdj \{v_1,v_2\}$,
		\item $(\Cn\G)^2_{\triangle}\cong\G^{2*}_{\triangle}$,
	\end{itemize}
	where $[v_1,v_2]$ denotes the closed edge between $v_1$ and $v_2$. We adopt this convention for the rest of the paper.
	
	To illustrate the usefulness of this approach, consider the following lemma.
	\begin{Lem}\label{trivsph}
		Let $\gamma\subset\G\setminus\{v_1,v_2\}$ be a cycle and suppose that there exist a path $p$ between $v_1$ and $v_2$ in $\G$ that does not intersect $\gamma$. Then a sphere $S^2$ that is a part of $S(\G\setminus\{v_1,v_2\})$ contributed by $\gamma$ maps homotopically trivially into $\G^{2*}_{\triangle}$. 
	\end{Lem}
	A path between vertices $u,v$ of a graph $\G$ is defined topologically as the subspace $\bigcup_{i=1}^n[v_i,v_{i+1}]$, where $[v_i,v_{i+1}]$, $1\leq i\leq n$ is a sequence of edges of $\G$, such that $v=v_1$ and $u=v_{n+1}$, and a \emph{cycle} is a path with $u=v$.
	\begin{proof}
		The path $p$ and cycle $\gamma$ are disjoint. Therefore $\gamma* p$ is a subspace of $\G^{2*}_{\triangle}$ homeomorphic to the $3$-dimensional ball $D^3$ (join of a circle and a line segment). A sphere contributed by $\gamma$ is given by $\gamma *\{v_1,v_2\}$ which is the boundary of $\gamma* p$. Thus, it is contractible inside $\G^{2*}_{\triangle}$.
	\end{proof}
	If every path $p$ between $v_1$ and $v_2$ intersects a cycle $\gamma$, then we say that $\gamma$ \emph{separates} $v_1$ and $v_2$.
	\section{Consequences of main theorem}
	
	\subsection{General facts} Recall that a graph is nondegenerate if it is connected and has a vertex of valence at least 3.
	\begin{Prop}\label{simconprop}
		Let $\G$ be a nondegenerate graph. Then, $(\Cn\G)^2_{\triangle}$ is simply connected. 
	\end{Prop}
	
	\begin{proof}
		By the assumption, there is a spanning tree $T$ of $\G$ that has a vertex of valence at least 3. We prove the claim by induction on the number of edges of $\G$ outside of $T$.
		
		From Corollary~\ref{corA} and Remark~\ref{treesrem}, we get:
		\[(\Cn T)^2_{\circ}\simeq S(T^2_{\circ})\simeq S\left(\bigvee_{i=1}^N S^1\right)\simeq \bigvee_{i=1}^N S^2,\]
		where
		\[N=\sum_{v\in V_T}(\nu(v)-1)(\nu(v)-2)-1.\]
		Since $T$ has a vertex of valence at least three, $N>1$.
		
		Consider a subgraph $\G'\subset\G$ containing $T$, and assume that $(\Cn\G')^2_{\triangle}$ is simply connected. We prove that $(\Cn(\ade(\G',v_1,v_2)))^2_{\triangle}$ is simply connected for any distinct $v_1,v_2\in V_{\G}=V_{\G'}$.
		
		By Corollary~\ref{maincor}, $(\Cn(\ade(\G',v_1,v_2)))^2_{\triangle}$ is a homotopy colimit of the diagram:
		\[*\leftarrow S(\G'\setminus\{v_1,v_2\})\rightarrow (\Cn\G')^2_{\triangle}\leftarrow S(\G'\setminus\{v_1,v_2\})\rightarrow *\]
		
		We first prove that the homotopy colimit of 
		\[*\leftarrow S(\G'\setminus\{v_1,v_2\})\rightarrow (\Cn\G')^2_{\triangle}\] 
		which we denote by $Z$, is simply connected and then analogous argument shows that 
		\[Z\leftarrow S(\G'\setminus\{v_1,v_2\})\rightarrow *\]
		is simply connected.
		
		We prove the first claim using the Van Kampen's theorem.
		We recall that the homotopy colimit of
		\[*\leftarrow S(\G'\setminus\{v_1,v_2\})\rightarrow (\Cn\G')^2_{\triangle}\]
		is given as a quotient of
		\[(S(\G'\setminus\{v_1,v_2\})\times[0,1])\sqcup (\Cn\G')^2_{\triangle},\]
		which squeezes down $S(\G'\setminus\{v_1,v_2\})\times\{0\}$ to a point and which identifies points $(x,1)\in S(\G'\setminus\{v_1,v_2\})\times\{1\}$ with the image of $x$ in $(\Cn\G')^2_{\triangle}$.
		
		Consider two open subsets $U,V\subset Z$. $U$ given as the image of \[\big(S(\G'\setminus\{v_1,v_2\})\times(0,1]\big)\sqcup (\Cn\G')^2_{\triangle}\] under the quotient map, and $V$ given as the image of $S(\G'\setminus\{v_1,v_2\})\times[0,1)$.
		\begin{itemize}
			\item $U\cap V$ is homotopy equivalent to $S(\G'\setminus\{v_1,v_2\})$, and hence it is path connected.
			\item $U$ is homotopy equivalent to $(\Cn\G')^2_{\triangle}$, and thus it is simply connected.
			\item $V$ is homotopy equivalent to a point, so it is also simply connected.
		\end{itemize}
		Van Kampen's theorem applies, and as a result, $Z$ is simply connected.
		
		We repeat the argument for 
		\[Z\leftarrow S(\G'\setminus\{v_1,v_2\})\rightarrow *\]
		to complete the proof.
	\end{proof}
	
	Notice that the maps from Corollary~\ref{maincor}:
	\begin{itemize}
		\item $\bigvee_{i=1}^{\nu(v)-1}S^1\to(\Cn\G)^2_{\triangle}$ and
		\item $\left(\bigvee_{i=1}^{b_0-1}S^1\right)\vee\left(\bigvee_{i=1}^{b_1}S^2\right) \rightarrow (\Cn\G)^2_{\triangle}$ restricted to $\bigvee_{i=1}^{b_0-1}S^1$
	\end{itemize}
	represent an element of $\pi_1((\Cn\G)^2_{\triangle})$, which for nondegenerate $\G$ is trivial as the proposition above states. Thus, we get:
	\begin{Cor}\label{attcor}
		For a nondegenerate $\G$:
		\begin{enumerate}[label=\roman*)]
			\item $(\Cn(\adv(\G,v)))^2_{\triangle}\simeq(\Cn\G)^2_{\triangle}\vee\bigvee_{i=1}^{2(\nu(v)-1)}S^2.$
			\item If $\left(\bigvee_{i=1}^{b_0-1}S^1\right)\vee\left(\bigvee_{i=1}^{b_1}S^2\right) \rightarrow (\Cn\G)^2_{\triangle}$ maps $\bigvee_{i=1}^{b_1}S^2$ to a contractible subspace, then:
			\[(\Cn(\ade(\G,v_1,v_2)))^2_{\triangle}\simeq(\Cn\G)^2_{\triangle}\vee\left(\bigvee_{i=1}^{2(b_0-1)}S^2\right)\vee\left(\bigvee_{i=1}^{2b_1}S^3\right) \]
		\end{enumerate}
		
	\end{Cor}
	\begin{proof}
		The homotopy colimit of the diagram in Corollary~\ref{maincor}~\ref{maincora} is obtained as a quotient of \[\Cn\left(\bigvee_{i=1}^{\nu(v)-1}S^1\right)\sqcup (\Cn\G)^2_{\triangle}\sqcup \Cn\left(\bigvee_{i=1}^{\nu(v)-1}S^1\right)\]
		by the relation that identifies the bases of the cones with images of respective maps $\bigvee_{i=1}^{\nu(v)-1}S^1\to(\Cn\G)^2_{\triangle}$  which are contractible inside $(\Cn\G)^2_{\triangle}$ as $\pi_1((\Cn\G)^2_{\triangle})=0$. Thus, the quotient is homotopy equivalent to $(\Cn\G)^2_{\triangle}\vee\bigvee_{i=1}^{2(\nu(v)-1)}S^2$. The second part is proved in the same way.
	\end{proof}
	\begin{Cor}
		If $\G\setminus\{v_1,v_2\}$ is contractible, then $\G^{2*}_{\triangle}\simeq(\ade(\G,v_1,v_2))^{2*}_{\triangle}$.
	\end{Cor}
	\begin{Prop}\label{discon}
		Let $\G_1,\dots,\G_n$ be connected components of $\G=\bigsqcup_{i=1}^n\G_i$ and suppose they are nondegenerate. Then $\G_{\triangle}^{2*}$ is homotopy equivalent to
		\[\bigvee_i(\G_i)^{2*}_{\triangle}\vee\bigvee_{i\neq j}\bigl((\G_i* \G_j)\vee S\G_i\vee S\G_j\bigr)\vee \bigvee_{k=1}^{(n-1)^2}S^1.\]
		In particular if Conjecture~\ref{conj} holds for all $\G_i$, then it holds for $\G$.  
	\end{Prop}
	\begin{proof}
		Let $X=\bigsqcup_{i=1}^nX_n$ be a disjoint union of simplicial complexes. Then,
		\[X^{2*}_{\triangle}=\bigcup_{1\leq i,j\leq n}X_i \sdj X_j.\]
		Notice that after identifying $X_i$ and $X_j$ with bases of the simplicial deleted join $X_i\sdj X_j$ we have:
		\[(X_i\sdj X_j) \cap (X_s\sdj X_t)=\begin{cases}
			(X_i)^{2*}_{\triangle} & \text{if }i=j=s=t\\
			X_i*X_j & \text{if } i=s\neq t=j\\
			X_i & \text{if } i=s \text{ and } j\neq t \\
			X_j & \text{if } i\neq s \text{ and } j=t \\
			\emptyset & \text{in all the other cases}
		\end{cases}\]
		The homotopy type of $X^{2*}_{\triangle}$ is obtained as a homotopy colimit of the diagram given by following spaces:
		\begin{itemize}
			\item two copies of each $X_i$ denoted $X_i^{(1)}$ and $X_i^{(2)}$ representing left and right side of the simplicial deleted join,
			\item a copy of $(X_i)^{2*}_{\triangle}$ for each $i$, and
			\item a copy of $X_i*X_j$ for each $i\neq j$.
		\end{itemize}
		and following maps, which are inclusions of subcomplexes: 
		\begin{itemize}
			\item $X^j_i\to (X_i)^{2*}_{\triangle}$, for $j=1,2$,
			\item $X^1_i\to X_i*X_j$,
			\item $X^2_i\to X_j*X_i$.
		\end{itemize}
		
		Let $X$ be a graph $\G$ and let $X_i=\G_i$ be as in the statement of the proposition. Notice that all maps in the diagram are nullhomotopic in this case as $\G_i$ are nondegenerate and Proposition~\ref{simconprop} applies. By Lemma~\ref{gmclem} in the appendix, we can replace maps in the diagram by homotopic constant maps without changing homotopy type of the homotopy colimit. We also replace each copy of $\G_i$ with a bouquet of circles representing its homotopy type, but we will still refer to it as $\G_i^{(\varepsilon)}$.
		Denote the homotopy colimit of the resulting diagram by $X$. It is obtained as the repeated wedge sum of following pieces:
		\begin{itemize}
			\item $\G_i^{(\varepsilon)}*\overline{n}$ for each $i\in\{1,\dots,n\}$ and each $\varepsilon\in\{1,2\}$, as each copy of $\G_i$ is a domain of $n$ maps,
			\item a copy of $(\G_i)^{2*}_{\triangle}$ for each $i$,
			\item and a copy of $\G_i*\G_j$ for each pair $i\neq j$,
		\end{itemize}
		
		In this space, we retract both copies of $\G_i$ seen as the bases of their joins with $\overline{n}$, mapping them onto their images in $(\G_i)^{2*}_{\triangle}$. This induces a deformation retraction of $X$.
		
		Next, we perform another deformation retraction on the resulting space that converge all the points in $\G_i\sdj \G_j$ that are images of all maps with it as a codomain to one chosen point.
		
		We rely on path-connectedness and CW-structure of each $\G_i\sdj \G_j$ to justify the existence of such deformations.
		
		The resulting space can be viewed as a graph, where vertices correspond to $\G_i\sdj \G_j$ and edges correspond to suspensions $S\G^{(\varepsilon)}_k$ that connect $\G_i\sdj \G_j$ to $(\G_k)^{2*}_{\triangle}$, where $i\neq j$ and $k$ is equal to either $i$ or $j$ (depending on $\varepsilon$). If a basepoint is chosen in each $\G^{(\varepsilon)}_k$, then a topological graph given by this graph structure can be viewed as a subspace of this space. Its vertices are the chosen points in $\G_i\sdj \G_j$ and edges are the suspensions of the chosen basepoints in $\G^{(\varepsilon)}_k$. Denote this topological graph by $G$.
		
		Contracting a spanning tree of $\G^{(\varepsilon)}_k$ induces a deformation retraction of $S\G^{(\varepsilon)}_k$ onto a bouquet of 2-spheres.
		
		We perform the final deformation retraction that deforms all $S\G^{(\varepsilon)}_k$'s into bouquets of $2$-spheres, moves them along the edge of $G$ to the image of $\G^{(\varepsilon)}_k$ in $(\G_k)^{2*}_{\triangle}$, and then contracts a spanning tree of $G$ to a point. $G$ has $n^2$ vertices and $2(n^2-n)$ edges, and therefore contracting its spanning tree results in the following bouquet of circles: $\bigvee_{k=1}^{(n-1)^2}S^1$.
		
		The resulting space is the wedge sum:
		
		\[\bigvee_i(\G_i)^{2*}_{\triangle}\vee\bigvee_{i\neq j}\bigl((\G_i* \G_j)\vee S\G_i\vee S\G_j\bigr)\vee \bigvee_{k=1}^{(n-1)^2}S^1.\]
	\end{proof}
	
	\subsection{Specific applications}
	\begin{Lem}\label{leafedge}
		Attaching a leaf or attaching an edge to vertices that already share one, in a graph satisfying the Conjecture~\ref{conj2}, yields a graph that also satisfies it.
	\end{Lem}
	\begin{proof}
		The part concerning a leaf follows from Corollary~\ref{attcor}. Suppose that $v_1$ and $v_2$ are connected by an edge in $\G$. No cycle $\gamma$ of $\G\setminus\{v_1,v_2\}$ separates $v_1$ and $v_2$ in $\G$ as the existing edge serves for the path $p$ for every $\gamma$. Thus, from Lemma~\ref{trivsph}, $\Cn\G\sdp\Cn\{v_1,v_2\}$ and $\Cn\{v_1,v_2\}\sdp\Cn\G$ are mapped homotopically trivially into $(\Cn\G)^2_{\triangle}$. As a result, Corollary~\ref{attcor} applies again, and
		\[(\Cn(\ade(\G,v_1,v_2)))^2_{\triangle}\simeq(\Cn\G)^2_{\triangle}\vee\left(\bigvee_{i=1}^{2(b_0-1)}S^2\right)\vee\left(\bigvee_{i=1}^{2b_1}S^3\right),\]
		where $b_0>1$.
	\end{proof}
	\begin{Rem}
		Trees with loops investigated in \cite{CheLut} can be obtained by applying these two operations. More precisely, a tree with loops is a graph formed as an iterated wedge sum of trees and circles. A wedge sum of a graph with a tree is obtained by repeatedly attaching leaves, while a wedge sum with a circle is realized by adding a leaf and subsequently doubling the newly added edge. 
		
		Consequently, Lemma~\ref{leafedge} implies that Conjecture~\ref{conj2} holds for all trees with loops, and the explicit homotopy type of their deleted squares can be computed recursively by iterating the wedge formula above. 
	\end{Rem}
	The following proposition is an analogue of Proposition~16 and Proposition~18 from~\cite{farber2010topology}.
	\begin{Prop}
		Suppose that nondegenerate $\G$ has a double edge, or a leaf. Then it is not mature, nor ripe.
	\end{Prop}
	\begin{proof}
		This is a consequence of the proof of Lemma~\ref{leafedge}, Proposition~\ref{propA}, and the fact that the zeroth Betti number of the complement of the set of two vertices joined by an edge is greater than 1.
	\end{proof}
	\begin{Prop}
		A graph containing two disjoint cycles is not unripe.
	\end{Prop}
	\begin{proof}
		Let $\G$ be a graph with disjoint cycles $\gamma_1$, $\gamma_2$. Let $[v_1,v_2]$ be an edge of $\gamma_1$ that does not lie in some spanning tree of $\G$. Let $\G'$ be obtained from $\G$ by removing (the interior of) $[v_1,v_2]$ so that $\G=\ade(\G',v_1,v_2)$. Let $p=\gamma_1\setminus[v_1,v_2]$. It is a path that assures that $\gamma_2$ does not separate $v_1$ and $v_2$. Therefore, by Lemma~\ref{trivsph} and Corollary~\ref{attcor}, $\gamma_2*\{v_1,v_2\}$ is mapped trivially into $(\G')^{2*}_{\triangle}$ and contributes a 3-sphere as a wedge sum component of $(\G)^{2*}_{\triangle}$. 
	\end{proof}
	The following proposition corresponds to Corollary~15 in \cite{farber2010topology}.
	\begin{Prop} 
		Let $\G$ be ripe and let $\hat{\G}=\G \cup e$ be obtained from $\G$ by adding an edge connecting two vertices $v_1,v_2\in\G$. If the complement $\G\setminus\{v_1,v_2\}$ is connected, then $\hat{\G}$ is ripe as well.
	\end{Prop}
	\begin{proof}
		If the complement $\G\setminus\{v_1,v_2\}$ is connected, then $S(\G\setminus\{v_1,v_2\})$ is homotopy equivalent to $\bigvee_{i=1}^{b_1(\G\setminus\{v_1,v_2\})}S^2$. Hence, from Corollary~\ref{maincor}~\ref{maincorb},
		$(\Cn\hat{\G})^2_{\triangle}$ is homotopy equivalent to \[\left(\bigvee_{i=1}^{2b_1(\G\setminus\{u,v\})}S^3\right)\vee(\Cn\G)^2_{\triangle}\]
		as $S(\G\setminus\{v_1,v_2\})$ maps nullhomotopically into $(\Cn\G)^2_{\triangle}$ due to ripeness of $\G$.
	\end{proof}
	The following proposition is analogous to Proposition~21 and Corollary~23 in \cite{farber2010topology}.
	\begin{Prop} 
		Suppose $n\geq 5$, then the complete graph $K_n$ is ripe. More precisely:
		
		\[(\Cn K_n)^2_{\triangle}\simeq \bigvee_{i=1}^{N+1}S^3,\]
		where \[N=\frac{1}{4}n(n-2)(n-3)(n-5).\]
	\end{Prop}
	\begin{proof}
		We prove the main claim by induction on $n$, using the general way of enlarging a complex to the case of $K_{n+1}$ obtained as
		$K_n\sqcup_V\Cn V$, where $V$ is the set of all vertices of $K_n$.
		
		Recall the diagram from Lemma~\ref{diaglem}:
		\[\begin{tikzcd}
			* & \Cn K_n\sdp\Cn V \arrow{l}\arrow{r} & (\Cn K_n)^2_{\triangle} & \Cn V\sdp\Cn K_n \arrow{l}\arrow{r} & *.
		\end{tikzcd}\]
		
		Homotopy type of $\Cn K_n\sdp\Cn V$ is obtained as a homotopy colimit of the diagram consisting of parts of the form
		\[\begin{tikzcd}
			\Cn(K_n\setminus\st(v)) & K_n\setminus\st(v) \arrow{l}\arrow{r} & K_n,
		\end{tikzcd}\]
		for each $v\in V$, where $\Cn(K_n\setminus\st(v))$ represents $\Cn K_n\sdp \{[v,0]\}$; $K_n\setminus\st(v)$ represents $\Cn K_n\sdp \quot{(\{v\}\times[0,1])}{\sim}$; and $K_n$ represents $\Cn K_n\sdp \{[v,1]\}$. In particular, all the parts share $K_n$.
		
		We can interpret this as coning in $K_n$ over complement of star of each vertex. We determine the homotopy type of this space.
		
		First, notice that $\Cn(K_n\setminus\st(v))\cong\Cn K_{n-1}$ has a deformation retraction that maps the apex of the cone to a chosen vertex of the base and edges connecting apex to the vertices of the base to edges emanating from the chosen vertex. There is no purely symmetric way to make this choice, but for computation purposes enumerate the vertices of $K_n$ with elements of $\mathbb{Z}_n$. Then for the deformation retraction of $\Cn(K_n\setminus\st(v_i))$ choose the vertex $v_{i+1}\in K_n\setminus\st(v_i)$.
		
		The retract is described as gluing in triangles into each subgraph of $K_n$ spanned by the set of vertices of the form $\{v_{i+1},v_j,v_k\}$, where $i$, $i+1$, $j$ and $k$ are distinct elements of $\mathbb{Z}_n$. Thus, the whole $\Cn K_n\sdp\Cn V$ is homotopy equivalent to a space obtained by gluing in triangles into $K_n$ along the subgraphs as above for every $i\in\mathbb{Z}_n$.
		
		This allows for counting over every such a subgraph. For a set $\{v_i,v_j,v_k\}$ of distinct vertices of $K_n$ there are one, two, or three triangles glued to it along the process:
		\begin{enumerate}[label=\roman*)]
			\item If $i,j,k$ are consecutive in cyclic order, then there is one triangle glued to $[\{v_i,v_j,v_k\}]$ and it comes from coning over $K_n\setminus\{\st(v_{i-1})\}$.
			There are exactly $n$ such sets $\{v_i,v_j,v_k\}$.
			\item If exactly two from $i,j,k$ are consecutive, without loss of generality, assume it is $j$ and $k$, then there are two triangles attached to $[\{v_i,v_j,v_k\}]$ that come from coning over complements of stars of $v_j$ and $v_k$.
			There are exactly $n(n-4)$ such sets $\{v_i,v_j,v_k\}$.
			\item The rest of the sets $\{v_i,v_j,v_k\}$ have pairwise nonconsecutive $i,j,k$, and therefore they receive three triangles.
		\end{enumerate}
		As every hollow triangle of $K_n$ is glued up at least once throughout this process, the resulting space is simply connected, and retracting its 1-skeleton to a point yields a bouquet of 2-spheres. Number of them is determined by the Euler characteristic.
		From a), b) and c), there are $3\binom{n}{3}-n(n-4)-2n$ 2-cells attached to $K_n$. Thus, the Euler characteristic of the whole space equals
		\[n-\binom{n}{2}+3\binom{n}{3}-n(n-4)-2n = \frac{1}{2}n(n-3)^2,\]
		and $\Cn K_n\sdp\Cn V$ is homotopy equivalent to a bouquet of $\frac{1}{2}n(n-3)^2-1$ spheres. We assume by induction that $(\Cn K_n)^2_{\triangle}$ has the homotopy type of a bouquet of spheres and thus has the trivial second homotopy group. In result:
		\[(\Cn K_{n+1})^2_{\triangle}\simeq (\Cn K_n)^2_{\triangle}\vee\bigvee_{i=1}^{n(n-3)^2-2}S^3.\]
		
		For the base case for induction, we show that $(K_5)^{2*}_{\triangle}$ is homeomorphic to $S^3$.
		
		$(K_5)^{2*}_{\triangle}$ can be identified with a subdivision of $(\Delta^4)^{(3)}$ (the 3-skeleton of the 4-simplex) constructed as follows.
		
		For new vertices of $(\Delta^4)^{(3)}$ take the barycenters of its 3-simplices. Denote a vertex that is a barycenter of 3-simplex opposite to a vertex $v$ of $\Delta^4$ by $v'$.
		
		Connecting the barycenter of a 3-simplex $\sigma$ with its vertices and with the barycenters of 3-simplices that share a 2-simplex with $\sigma$ yields a 1-skeleton of this subdivision.
		
		To complete the construction, take the flag completion of the obtained 1-complex.
		
		Isomorphism with $(K_5)^{2*}_{\triangle}$ is given on vertices. Vertex $v$ of $\Delta^4$ is identified with $(v,-1)$ in $(K_5)^{2*}_{\triangle}$ and new vertex $v'$ is identified with $(v,1)$.
		
		Hence, for $n>5$:
		
		\[(\Cn K_n)^2_{\triangle}\simeq \bigvee_{i=1}^{N+1}S^3,\]
		where \[N=\sum_{k=5}^{n-1}(k(k-3)^2-2)=\frac{1}{4}n(n-2)(n-3)(n-5).\]
	\end{proof}
	\section{Appendix}
	\subsection{Proof of Lemma~\ref{lemret}}
	Let $X$ be a simplicial complex with $n+1$ vertices. We consider it as a subcomplex of the standard simplex $\Delta^n$ in $\mathbb{R}^{n+1}$. For a pair of points $(x,y)\in X^2_{\triangle}$ with barycentric coordinates $x=(t_1,\dots,t_{n+1})$ and $y=(s_1,\dots,s_{n+1})$, consider hyperplanes determined by fixing $\min \{t_i,s_i\}$ on $i$-th coordinate. Inside $\Delta^n$ they bound an $n$-simplex, which we denote $\sigma(x,y)$.
	
	Notice that $\sigma(x,y)$ is the smallest simplex that contains both $x$ and $y$ and is similar to $\Delta^n$ through a homothety of positive ratio. Let $T_{x,y}\colon\mathbb{R}^{n+1}\to\mathbb{R}^{n+1}$ be this homothety, and let $k\geq1$ be its ratio, assuming that $T_{x,y}$ maps $\sigma(x,y)$ to $\Delta^n$. Then, the deformation retraction $f\colon\Conf_2(X)\to X^2_{\triangle}$ is given by \[f((x,y))=(T_{x,y}(x),T_{x,y}(y)).\]
	The homotopy with identity is given by $f_t((x,y))=(T^t_{x,y}(x),T^t_{x,y}(y))$, where $T^t_{x,y}$ is the homothety with the same center as $T_{x,y}$ and with the ratio $1+t(k-1)$.
	
	To see that $X^2_{\triangle}$ is a strong deformation retract, notice that in $\sigma(x,y)$, thought of as simplicial complex, $\supp(x)\cap\supp(y)=\emptyset$. Otherwise there would be a vertex $v$ of $\sigma(x,y)$ inside $\supp(x)\cap\supp(y)$ and that would mean that neither $x$ or $y$ lie in the opposite to $v$ face of $\sigma(x,y)$, which is impossible, because codimension one faces of $\sigma(x,y)$ lie in hyperplanes passing through either $x$ or $y$. Therefore $\supp(T_{x,y}(x))\cap\supp(T_{x,y}(y))=\emptyset$ in $\Delta^n$, and hence in $X$.
	
	On the other hand, if $(x,y)\in X^2_{\triangle}$, then $\sigma(x,y)=\Delta^n$ as $\min\{t_i,s_i\}=0$ for all $i$, because $\supp(x)$ and $\supp(y)$ in $X$ share no vertex of $X$ and thus no vertex of $\Delta^n$, and in that case $T_{x,y}$ is the identity map.
	
	The case of a cube complex can be dealt with similarly, but it requires more tedious technical details as for a pair of points in a cube, the analogous family of hyperplanes would not necessarily bound a cube, and thus there is no obvious map to build a deformation retraction on. For our purpose the case of simplicial complex is sufficient.
	$\square$
	\subsection{Proof of Lemma~\ref{simhoco}}
	First we prove a generalization of a classic fact on mapping cones (see \cite[page~88]{arkowitz2011introduction}).
	\begin{Lem}\label{gmclem}
		Let $X_{\G}$ be a diagram of spaces and let $a$ be an edge of $\G$. Consider a new diagram $Y_{\G}=(Y_v,\psi_e)_{\G}$ over the same graph $\G$ with spaces $Y_v=X_v$ for each vertex $v$, and maps $\psi_e=\varphi_e$ for each edge $e\neq a$, and let $\psi_a$ be a map homotopic to $\varphi_a$ by $H\colon X_{0_a}\times[0,1]\to X_{1_a}$, where $H(x,0)=\varphi_a(x)$ and $H(x,1)=\psi_a(x)$. Then $\gmc(X_{\G})$ and $\gmc(Y_{\G})$ are homotopy equivalent. 
	\end{Lem}
	\begin{proof}
		Define $H'\colon X_{0_a}\times[0,1]\to X_{1_a}$ by $H'(x,t)=H(x,1-t)$.
		
		Denote by $[x]_X$ and $[x,t_e]_X$ the equivalence classes of $x$ and $(x,t_e)$ in $\gmc(X_{\G})$ and analogously $[x]_Y$ and $[x,t_e]_Y$ in $\gmc(Y_{\G})$.
		
		Define maps $F\colon \gmc(X_{\G})\to \gmc(Y_{\G})$ and $F'\colon \gmc(Y_{\G})\to \gmc(X_{\G})$ by
		\begin{itemize}
			\item $F([x]_X)=[x]_Y$ and $F'([x]_Y)=[x]_X$ for each $v$ and each $x\in X_v=Y_v$
			\item $F([x,t_e]_X)=[x,t_e]_Y$ and $F'([x,t_e]_Y)=[x,t_e]_X$ for each $e\neq a$, each $x\in X_{0_e}$ and each $t\in[0,1]$
			\item $F([x,t_a]_X)=\begin{cases}
				[x,(2t)_a]_Y & \text{for } t\in[0,\frac{1}{2}]\\
				[H'(x,2t-1)]_Y & \text{for } t\in[\frac{1}{2},1]
			\end{cases}$
			\item $F'([x,t_a]_Y)=\begin{cases}
				[x,(2t)_a]_X & \text{for } t\in[0,\frac{1}{2}]\\
				[H(x,2t-1)]_X & \text{for } t\in[\frac{1}{2},1]
			\end{cases}$
		\end{itemize}
		Those maps are properly defined. Indeed,
		\begin{itemize}
			\item $F([x,(\frac{1}{2})_a]_X)=[x,1_a]_Y=[\psi_a(x)]_Y=[H'(x,0)]_Y$.
			\item $F([x,1_e]_X)=[x,1_e]_Y=[\psi_e(x)]_Y=[\varphi_e(x)]_Y=F([\varphi_e(x)]_X)$ for $e\neq a$ and $F([x,1_a]_X)=[\varphi_a(x)]_Y=F([\varphi_a(x)]_X)$.
		\end{itemize}
		Analogous argument works for $F'$.
		
		We show that $F'F$ is homotopic to $\id_{\gmc(X_{\G})}$, and analogous argument for $FF'$ proves that $F$ and $F'$ are homotopy equivalences.
		
		First, observe that:
		\begin{itemize}
			\item $F'F([x]_X)=[x]_X$ for each $v$ and each $x\in X_v$.
			\item $F'F([x,t_e]_X)=[x,t_e]_X$ for each $e\neq a$, each $x\in X_{0_e}$ and each $t\in[0,1]$.
			\item $F'F([x,t_a]_X)=\begin{cases}
				[x,(4t)_a]_X, & \text{for } t\in[0,\frac{1}{4}],\\
				[H(x,4t-1)]_X, & \text{for } t\in[\frac{1}{4},\frac{1}{2}],\\
				[H'(x,2t-1)]_X, & \text{for } t\in[\frac{1}{2},1].
			\end{cases}$
		\end{itemize}
		Define the homotopy $\mathcal{H}\colon \gmc(X_{\G})\times[0,1]\to \gmc(X_{\G})$ by
		\begin{itemize}
			\item $\mathcal{H}([x]_X,s)=[x]_X$ for each $v$, each $x\in X_v$, and each $s\in[0,1]$.
			\item $\mathcal{H}([x,t_e]_X,s)=[x,t_e]_X$ for each $e\neq a$, each $x\in X_{0_e}$ each $t\in[0,1]$, and each $s\in[0,1]$.
			\item $\mathcal{H}([x,t_a]_X,s)=\begin{cases}
				[x,(\frac{4t}{1+3s})_a]_X, & \text{for } t\in[0,\frac{1+3s}{4}],\\
				[H(x,4t-3s-1)]_X, & \text{for } t\in[\frac{1+3s}{4},\frac{1+s}{2}],\\
				[H'(x,2t-1)]_X, & \text{for } t\in[\frac{1+s}{2},1].
			\end{cases}$
		\end{itemize}
		Notice that $\mathcal{H}_0=F'F$ and $\mathcal{H}_1=\id_{GMC(X_{\G})}$.
		So this is the desired homotopy.
		
		To be sure, check that $4t-3s-1$ for $t\in[\frac{1+3s}{4},\frac{1+s}{2}]$ sits in $[0,1-s]$, and $2t-1$ for $t\in[\frac{1+s}{2},1]$ is in $[s,1]$. Also for $s=1$ and $t=1$, we have: \[[x,1_a]_X=[\varphi_a(x)]_X=[H(x,0)]_X=[H'(x,1)]_X.\]
	\end{proof}
	\begin{Cor}\label{gmccor}
		If two diagrams differ by homotopic maps, their generalized mapping cylinders are homotopy equivalent.
	\end{Cor}
	We also require a fundamental property of homotopy colimits, for the proof of which we refer to \cite[page~265]{kozlov2008combinatorial}.
	\begin{Lem}\label{hocolem}
		If a commutative map of diagrams of spaces consists of homotopy equivalences, then homotopy colimits of those diagrams are homotopy equivalent.
	\end{Lem}
	With those two facts we can prove Lemma~\ref{simhoco}.
	\begin{proof}
		Let $f\colon X_{\G}\to Y_{\G}$ be a homotopy commutative map of source--drain diagrams $X_{\G}=(X_v,\varphi_e)_{\G}$ and $Y_{\G}=(Y_v,\psi_e)_\G$, such that all $f_v$ are homotopy equivalences.
		
		Define diagrams $Z_{\G}=(Z_v,\alpha_e)_{\G}$ and $Z'_{\G}=(Z'_v,\beta_e)_{\G}$ by setting $Z_d=Z'_d=Y_d$ for all drains $d$, $Z_s=Z's=X_s$ for all sources $s$, and $\alpha_e=f_{1_e}\varphi_e$, $\beta_e=\psi_ef_{0_e}$ for all edges $e$.
		
		Define maps of diagrams $g\colon X_{\G}\to Z_{\G}$ and $g'\colon Z'_{\G}\to Y_{\G}$ by $g_d=f_d$ and $g'_d=\id_{Y_d}$ for all drains $d$ and by $g_s=\id_{X_s}$ and $g'_s=f_s$ for all sources $s$.
		
		Notice that they are both commutative. Indeed, for each edge $e$ the vertex $0_e$ is a source and $1_e$ is a drain from the source--drain property. Therefore the squares
		\[\begin{tikzcd}
			X_{0_e}\arrow{r}{\varphi_e}\arrow{d}{g_{0_e}} & X_{1_e}\arrow{d}{g_{1_e}} \\
			Z_{0_e}\arrow{r}{\alpha_e} & Z_{1_e}
		\end{tikzcd}\begin{tikzcd}
			Z'_{0_e}\arrow{r}{\beta_e}\arrow{d}{g_{0_e}} & Z'_{1_e}\arrow{d}{g_{1_e}} \\
			Y_{0_e}\arrow{r}{\psi_e} & Y_{1_e}
		\end{tikzcd}\]
		become
		\[\begin{tikzcd}
			X_{0_e}\arrow{r}{\varphi_e}\arrow{d}{\id} & X_{1_e}\arrow{d}{f_{1_e}} \\
			X_{0_e}\arrow{r}{f_{1_e}\varphi_e} & Y_{1_e}
		\end{tikzcd}\begin{tikzcd}
			X_{0_e}\arrow{r}{\psi_ef_{0_e}}\arrow{d}{f_{0_e}} & Y_{1_e}\arrow{d}{\id} \\
			Y_{0_e}\arrow{r}{\psi_e} & Y_{1_e}
		\end{tikzcd}\]
		and they obviously commute.
		By the assumption, $\alpha_e$ is homotopic to $\beta_e$ for all $e$ and $Z_v=Z'_v$ for all $v$. Thus, from Corollary~\ref{gmccor}, homotopy colimits of $Z_{\G}$ and $Z'_{\G}$ are homotopy equivalent.
		
		By Lemma~\ref{hocolem}, the homotopy colimit of $X_{\G}$ is homotopy equivalent to the homotopy colimit of $Z_{\G}$ and the homotopy colimit of $Z'_{\G}$ is homotopy equivalent to the homotopy colimit of $Y_{\G}$.
		
		Therefore homotopy colimits of $X_{\G}$ and $Y_{\G}$ are homotopy equivalent.
	\end{proof}
	\subsection{Homotopy type of configuration space of two particles in a tree}
	Homotopy type of the configuration space of $n$ particles in a tree was studied in \cite{Ghrist}. In \cite{farber2005collision} an explicit formula for the case of 2 particles was given and it was fully proven in \cite{FarberConf}. We show an elementary proof for the convenience of the reader.
	
	Let $T$ be a tree. Denote by $N_2(T)$ the subcomplex of $T^2$ consisting of all closed cells intersecting the diagonal. Consider the subcomplex $N_2(T) \cap T^2_{\triangle}$ of $T^2$. This subcomplex is a graph. Its vertices are ordered pairs of distinct vertices of $T$ lying in the same edge or in two edges sharing a vertex (which is not in the pair), and its edges are ordered pairs consisting of an edge $e$ and vertex $v$ of $T$ that is incident to an edge sharing exactly one vertex with $e$. Denote this graph by $\partial T^2_{\triangle}$. It is isomorphic to the barycentric subdivision of a bipartite graph with set of vertices $E_T \sqcup E_T$ and edges between vertices corresponding to edges sharing a vertex in $T$.
	
	We show that $\partial T^2_{\triangle}$ is a deformation retract of $T^2_{\triangle}$. Consider a natural metric $d$ on $T$ such that each edge is isometric to the unit interval $[0,1]$. With such metric, $T$ is a uniquely geodesic space. That means that for any two points $x, y\in T$, there is exactly one path between $x$ and $y$ of shortest length called a geodesic. In particular, there is an isometric embedding $I\colon [0,d(x,y)]\to T$ with $I(0)=x$ and $I(d(x,y))=y$. Let $m=I\left(\frac{d(x,y)}{2}\right)$ be the midpoint of this path. Let $v_x\in I\left(\left(0,\frac{d(x,y)}{2}\right)\right)$ and $v_y\in I\left(\left(\frac{d(x,y)}{2},d(x,y)\right)\right)$ be the two vertices on the geodesic between $x$ and $y$ that are closest to $m$. Let $M=\max\{d(v_x,m),d(v_y,m)\}$ and put
	\[f((x,y))=\left(I\left(\frac{d(x,y)}{2}-M\right), I\left(\frac{d(x,y)}{2}+M\right)\right).\]
	The map $f\colon T^2_{\triangle}\to T^2_{\triangle}$ is the desired deformation retraction.
	
	We determine the homotopy type of $\partial T^2_{\triangle}$ using the Euler characteristic of this space:
	\[\chi(\partial T^2_{\triangle})=2\#E_T - 2\sum_{v\in V_T} \binom{\nu(v)}{2} = \sum_ {v\in V_T} \nu(v)(2-\nu(v)).\]
	We have proved the following
	\begin{Lem}\label{trees}
		$T^2_{\triangle}$ has the homotopy type of a bouquet of $1 + \sum_ {v\in V_T} \nu(v)(\nu(v)-2)$ circles.
	\end{Lem}
	\begin{Rem}\label{treesrem}
		Since $\#E_T=\#V_T-1$:
		\begin{align*}
			1 + \sum_{v\in V_T} \nu(v)(\nu(v)-2) 
			&= 1 + \sum_{v\in V_T} \nu(v)^2 - 4\#E_T \\
			&= \sum_{v\in V_T} \nu(v)^2 - 6\#E_T + 2\#V_T - 1 \\
			&= \sum_{v\in V_T}(\nu(v)-1)(\nu(v)-2) - 1.
		\end{align*}
		Therefore our formula agrees with the one in \cite{farber2005collision}.
	\end{Rem}
	\section{Acknowledgements}
	This work forms part of the author's Ph.D. research at the University of Wroc\l{}aw.
	
	The author is deeply grateful to his advisor, \'Swiatos\l{}aw Gal, for his help and guidance, and in particular for proposing the conjecture studied here.
	
	The author is also grateful to Ben Knudsen for insightful comments and questions regarding an earlier draft of this work.
	
	\bibliographystyle{plain}

\end{document}